\documentclass[nosumlimits]{amsart}
\usepackage{amscd,amssymb,epsfig}
\theoremstyle{plain}
\newtheorem{prop}[subsubsection]{Proposition}
\newtheorem{thm}[subsubsection]{Theorem}
\newtheorem{lem}[subsubsection]{Lemma}
\newtheorem{cor}[subsubsection]{Corollary}
\newtheorem*{thm*}{Theorem}

\theoremstyle{remark}
\newtheorem{rem}[subsubsection]{Remark}
\newtheorem*{rem*}{Remark}
\newtheorem*{ack}{Acknowledgments}

\newtheorem*{conv}{Conventions}

\theoremstyle{definition}
\newtheorem{defn}[subsubsection]{Definition}
\newtheorem{exm}[subsubsection]{Example}

\numberwithin{equation}{section}

\newcommand{\A}{{\mathcal A}}
\newcommand{\B}{{\mathcal B}}
\newcommand{\LL}{{\mathcal L}}
\newcommand{\W}{{\mathcal W}}
\newcommand{\Z}{{\mathbb Z}}
\newcommand{\R}{{\mathbb R}}
\newcommand{\C}{{\mathbb C}}
\newcommand{\Q}{{\mathbb Q}}
\newcommand{\F}{{\mathbb F}}
\newcommand{\ii}{{\sqrt{-1}}}

\renewcommand{\a}{{\alpha }}
\renewcommand{\b}{{\beta }}
\renewcommand{\c}{{\gamma }}
\newcommand{\bm}{{\mathbf m}}
\newcommand{\bx}{{\mathbf x}}
\newcommand{\by}{{\mathbf y}}
\newcommand{\z}{{\zeta }}
\newcommand{\s}{{\sigma }}
\renewcommand{\ll}{{\ell }}
\renewcommand{\H}{{\mathcal H}}

\DeclareMathOperator{\rank}{rank}
\DeclareMathOperator{\conj}{conj}
\DeclareMathOperator{\id}{id}
\DeclareMathOperator{\Aut}{Aut}
\DeclareMathOperator{\Out}{Out}
\DeclareMathOperator{\BB}{B}

\DeclareMathOperator{\codim}{codim}
\DeclareMathOperator{\Real}{Re}
\DeclareMathOperator{\Imag}{Im}
\DeclareMathOperator{\gr}{gr}

\begin{document}

\title[fiber-type arrangements \& orbit configuration spaces]
{Monodromy of Fiber-Type Arrangements and Orbit Configuration Spaces}
\author{Daniel C.~Cohen}
\address{Department of Mathematics, Louisiana State University, 
Baton Rouge, LA 70803}
\email{cohen@math.lsu.edu}
\urladdr{http://math.lsu.edu/\~{}cohen}
\thanks{Partially supported by grants LEQSF(1996-99)-RD-A-04 and 
LEQSF(1999-2002)-RD-A-01 from the Louisiana Board of Regents.}



\begin{abstract}
We prove similar theorems concerning the structure of bundles 
involving complements of fiber-type hyperplane arrangements and orbit 
configuration spaces.  These results facilitate analysis of the 
fundamental groups of these spaces, which may be viewed as 
generalizations of the Artin pure braid group.  In particular, we 
resolve two disparate conjectures.  We show that the Whitehead group 
of the fundamental group of the complement of a fiber-type arrangement 
is trivial, as conjectured by Aravinda, Farrell, and Rouchon 
\cite{AFR}.  For the orbit configuration space corresponding to the 
natural action of a finite cyclic group on the punctured plane, we 
determine the structure of the Lie algebra associated to the lower 
central series of the fundamental group.  Our results show that this 
Lie algebra is isomorphic to the module of primitives in the homology 
of the loop space of a related orbit configuration space, as 
conjectured by Xicot\'encatl~\cite{Xi}.
\end{abstract}

\maketitle

\section*{Introduction} \label{sec:intro}

Let $M$ be a manifold without boundary of dimension at least two.  The 
{\em configuration space} of $n$ ordered points in $M$ is the subspace 
of the product space $M^{n}$ defined by
\begin{equation*} \label{eq:ConfigSpaceM}
F(M,n) = \{(x_{1},\dots,x_{n})\in M^{n} \mid x_{i} \neq x_{j}\ 
\text{if}\ i \neq j\}.
\end{equation*}
These spaces arise in numerous contexts, including of course that of 
braid groups.  The symmetric group $\Sigma_n$ acts freely on $F(M,n)$.  
The fundamental group of $F(M,n)/\Sigma_n$ is called the (full) braid 
group of $M$, and that of $F(M,n)$ is the pure braid group of $M$.  In 
the case $M=\C$, these groups are the classical Artin braid groups, 
see the books of Birman \cite{Bi} and Hansen \cite{Ha}, to which we 
refer as general references on braids.
 
Let $Q_n$ denote a set of $n$ distinct points in $M$.  A fundamental 
property of configuration spaces is given by the following classical 
result, which we will use extensively.
\begin{thm*}[Fadell and Neuwirth \cite{FN}, Theorem~3]
\label{thm:FN}
For $\ll \le n$, the projection onto the first $\ll$ coordinates, 
$p:F(M,n)\to F(M,\ll)$, is a locally trivial bundle, with fiber 
$F(M\setminus~{\hskip -3pt Q_\ll},n-\ll)$.
\end{thm*}

The focus of this paper is on two related generalizations of 
configuration spaces.  First, we study the class of {\em fiber-type 
arrangements}.  An arrangement of hyperplanes is a finite collection 
of codimension one affine subspaces of Euclidean space $\C^{\ll}$.  
See Orlik and Terao \cite{OT} as a general reference on arrangements.  
The complement of an arrangement $\A$ is the manifold 
$M(\A)=\C^{\ll}\setminus\bigcup_{H\in\A}H$.  The configuration space 
$F(\C,n)$ may be realized as the complement in $\C^{n}$ of the braid 
arrangement consisting of the hyperplanes $H_{i,j}=\ker(x_{j}-x_{i})$, 
$1\le i<j \le n$.  Briefly, an arrangement is fiber-type if, as is the 
case for the braid arrangement, the complement sits atop a tower of 
fiber bundles, the projection maps of which are the restrictions to 
hyperplane complements of linear maps $\C^{k}\to\C^{k-1}$.

The second generalization arises in the following way.  Let $\Gamma$ 
be a (finite) group, acting freely on the manifold $M$.  The {\em 
orbit configuration space} $F_\Gamma(M,n)$ is the subspace of the 
product space $M^n$ consisting of all ordered $n$-tuples 
$(x_1,\dots,x_n)$ of points in $M$ for which the orbits $\Gamma\cdot 
x_i$ and $\Gamma\cdot x_j$ do not intersect for $i\neq j$.  These 
spaces were recently studied by Xicot\'encatl \cite{Xi}, who showed 
that, like their classical counterparts, they support certain Lie 
algebra structures in their loop space homology, and proved a 
fibration theorem generalizing the Fadell-Neuwirth theorem stated 
above.

In this paper, we use the Fadell-Neuwirth theorem to determine the 
structure of certain bundles involving both complements of fiber-type 
arrangements and orbit configuration spaces.  These results, the 
proofs of which are straightforward and virtually identical, 
facilitate analysis of the fundamental groups of these spaces.  As 
these spaces generalize configuration spaces, their fundamental groups 
may be viewed as generalizations of the Artin pure braid group.  We 
pursue several structural aspects of these generalized pure braid 
groups.  Much is already known about a number of these groups.

For a fiber-type arrangement $\A$, the group $G=\pi_1(M(\A))$ is an 
``almost direct product'' of free groups.  That is, the group $G$ may 
be realized as an iterated semidirect product of free groups, the 
factors of which act on one another by conjugation.  In \cite{FR1}, 
Falk and Randell use this structure to prove that fiber-type 
arrangement groups satisfy the celebrated LCS formula, relating the 
Betti numbers and the ranks of the lower central series quotients of 
these groups.  For the Artin pure braid group, this result was 
obtained by other means by Kohno \cite{Ko}.  Using the almost direct 
product structure of the group $G$ of a fiber-type arrangement, one 
can also construct a finite free resolution of the integers over the 
group ring ${\Z}G$ and show that $G$ is of type FL \cite{CScc}, show 
that $G$ is residually nilpotent \cite{FR2}, orderable \cite{Pa}, 
etc..

Our results further reveal the structure of these groups.  From the 
relationship with configuration spaces we establish, it follows that 
the conjugation action in the almost direct product structure is given 
by pure braid automorphisms.  Using this, we show that fiber-type 
arrangement groups are strongly poly-free, see Definition 
\ref{def:SPF}.  Aravinda, Farrell, and Roushon \cite{AFR} have 
recently shown that the Whitehead group of any strongly poly-free 
group is trivial.  Thus, the Whitehead group of a fiber-type 
arrangement group is trivial, as conjectured in \cite{AFR}.

In principle, we give an algorithm for presenting the group of a 
fiber-type arrangement as an almost direct product.  Using a natural 
generalization of the techniques for computing the braid monodromy of 
a complex line arrangement developed in \cite{CSbm}, we provide a 
method for calculating the pure braids which dictate the structure of 
fiber-type arrangement groups.  We illustrate the method using the 
Coxeter arrangement of type $\BB$, giving a presentation of the 
Brieskorn generalized pure braid group $PB_n$ which exhibits the 
iterated semidirect product structure of this group.

The complex reflection arrangements associated to full monomial groups 
provide a bridge between the two generalizations of configuration 
spaces we consider.  The complement of such an arrangement may be 
realized as the orbit configuration space $F_\Gamma(\C^*,n)$, where 
$\Gamma$ is a finite cyclic group, acting on $\C^*$ by multiplication 
by a primitive root of unity.  By combining the above methods with 
known results on fundamental groups of hyperplane complements, we also 
obtain a presentation for the fundamental group of this orbit 
configuration space.  This presentation is used to study the Lie 
algebra associated to the lower central series of this ``pure monomial 
braid group,'' revealing that the structure of this Lie algebra is as 
conjectured by Xicot\'encatl \cite{Xi}.

\begin{conv} Denote by $\Aut(G)$ the group of right automorphisms of a 
group $G$, with multiplication $\a\cdot\b=\b\circ\a$.  For $u,v\in G$, 
write $u^v=v^{-1}uv$ and $[u,v]=uvu^{-1}v^{-1}$.
\end{conv}

\section{Fiber-Type Arrangements} \label{sec:fta}

\subsection{Monodromy of Fiber-Type Arrangements} 
\label{subsec:ftmono}
\ 
\medskip

In this section, we identify the monodromy of a strictly linearly 
fibered arrangement.  When iteratively applied to fiber-type 
arrangements, this identification gives the iterated semidirect 
product structure of the fundamental group of the complement of such 
an arrangement.  We first recall the definitions of these 
arrangements, see \cite{FR1,OT}.

\begin{defn} \label{defn:slfdef}
A hyperplane arrangement $\A$ in $\C^{\ll+1}$ is {\em strictly 
linearly fibered} if there is a choice of coordinates 
$(\bx,z)=(x_1,\dots,x_\ll,z)$ on $\C^{\ll+1}$ so that the restriction, 
$p$, of the projection $\C^{\ll+1}\to\C^\ll$, $(\bx,z)\mapsto \bx$, to 
the complement $M(\A)$ is a fiber bundle projection, with base 
$p(M(\A))=M(\B)$, the complement of an arrangement $\B$ in $\C^\ll$, 
and fiber the complement of finitely many points in $\C$.  We say $\A$ 
is strictly linearly fibered over $\B$.
\end{defn}

\begin{defn} \label{defn:ftdef}
An arrangement $\A=\A_1$ of finitely many points in $\C^1$ is {\em 
fiber-type}.  An arrangement $\A=\A_\ll$ of hyperplanes in $\C^\ll$ is 
{\em fiber-type} if $\A$ is strictly linearly fibered over a 
fiber-type arrangement $\A_{\ll-1}$ in $\C^{\ll-1}$.
\end{defn}

The complement of a fiber-type arrangement sits atop a tower of fiber 
bundles 
\begin{equation*} \label{eq:BundleTower}
M(\A_\ll) \xrightarrow{p_{\ll}} M(\A_{\ll-1}) 
\xrightarrow{p_{\ll-1}} \cdots 
\xrightarrow{p_{2}} M(\A_1) = \C \setminus \{d_1
\text{ points}\},
\end{equation*}
where the fiber of $p_{k}$ is homeomorphic to the complement of 
$d_{k}$ points in $\C$.  Repeated application of the homotopy exact 
sequence of a bundle shows that $M(\A_{\ll})$ is a $K(\pi,1)$ space, 
where $\pi=\pi_{1}(M(\A_{\ll}))$.  The integers $\{d_1,\dots,d_\ll\}$ 
are called the {\em exponents} of the fiber-type arrangement $\A$.  In 
general, an arrangement $\A$ is said to be $K(\pi,1)$ if the 
complement $M(\A)$ is an Eilenberg-MacLane space of type $K(\pi,1)$.

\begin{rem} \label{rem:supersolvable}
The rank of an arrangement $\A$ is the largest number of linearly 
independent hyperplanes in $\A$.  A useful alternative definition of a 
fiber-type, or supersolvable, arrangement is given in \cite{FT}: An 
(affine) arrangement $\A$ is {\em supersolvable} if there is a 
sequence $\A=\A_\ll \supseteq \cdots \supseteq \A_1$ such that 
$\rank\A_{j+1}=\rank\A_{j}+1$ for each $j$, and for distinct $H,H' \in 
\A_j$ with $H\cap H' \neq \emptyset$, there exists $H'' \in \A_i$ with 
$i<j$ and $H\cap H' \subset H''$.
\end{rem}

\begin{exm} \label{exm:braidA}
Consider the braid arrangement, $\{\ker(y_i-y_j), 1\le i<j\le n\}$, in 
$\C^n$, with complement $F(\C,n)=\{\by\in\C^n \mid y_i \neq y_j\ 
\text{if}\ i\neq j\}$, the configuration space of $n$ ordered points 
in $\C$, where $\by=(y_1,\dots,y_n)$.  The braid arrangement is the 
prototypical example of a fiber-type arrangement.  By the 
Fadell-Neuwirth theorem, projection onto the first $n$ coordinates 
yields a bundle $p_{n+1}:F(\C,n+1)\to F(\C,n)$, with fiber the 
complement of $n$ points in $\C$.

As noted previously, the fundamental groups of $F(\C,n)$ and 
$F(\C,n)/\Sigma_n$ are the classical Artin braid groups $P_n$ and 
$B_n$, see \cite{Bi,Ha}.  We record presentations of these groups.  
The full braid group has presentation
\begin{equation} \label{eq:FullBraidPres}
B_n=\Biggl\langle \s_i\ (1\le i < n) \ \Biggm|
\begin{split}
\s_i \s_{i+1} \s_i &= \s_{i+1} \s_i \s_{i+1}\ (1\le i < n-1)\\
\s_i \s_j &= \s_j \s_i\ (|j-i|>1) 
\end{split}\ 
\Biggr\rangle.
\end{equation}

The bundle of configuration spaces $p_{n+1}:F(\C,n+1)\to F(\C,n)$ 
admits a section $s:F(\C,n)\to F(\C,n+1)$, given by $s(\by)=(\by,z)$, 
where $z=\bigl(1+\sum_{i=1}^n|y_i|^2\bigr)^{1/2}$.  From this and the 
homotopy sequence of the bundle, it follows that the pure braid group 
admits the structure of an iterated semidirect product of free groups, 
$P_n=\F_{n-1} \rtimes_{\a_{n-1}} \dots \rtimes_{\a_2} \F_1$.  The 
monodromy homomorphism $\a_k:P_k\to\Aut(\F_k)$ is the restriction to 
$P_k$ of the Artin representation $\a_k:B_k\to\Aut(\F_k)$, defined by
\begin{equation} \label{eq:ArtinRep}
\a_k(\s_i)(t_j)=\begin{cases}
t_i^{} t_{i+1}^{} t_i^{-1} & \text{if $j=i$,}\\
t_i^{} & \text{if $j=i+1$,}\\
t_j^{} & \text{otherwise,}
\end{cases}
\end{equation}
where $\F_k=\langle t_1,\dots,t_k\rangle$.  The standard presentation 
of the pure braid group exhibits the iterated semidirect product 
structure.  The group $P_n$ has generators
\begin{equation} \label{eq:PureBraidGens}
A_{i,j}^{}=\s_{j-1}^{}\cdots\s_{i+1}^{} 
\s_i^2 \s_{i+1}^{-1}\cdots\s_{j-1}^{-1} \quad (1\le i<j\le n),
\end{equation}
and, for $j<l$, defining relations
\begin{equation} \label{eq:PureBraidRels}
A_{i,j}^{-1}A_{k,l}^{}A_{i,j}^{}=
\begin{cases}
(A_{i,l}A_{j,l})A_{k,l}(A_{i,l}A_{j,l})^{-1}&
\text{if $k=i$ or $k=j$,}\\
{[}A_{i,l}A_{j,l}{]} A_{k,l} {[}A_{i,l}A_{j,l}{]}^{-1}&
\text{if $i<k<j$,}\\
A_{k,l}&\text{otherwise.}
\end{cases}
\end{equation}
\end{exm}

We now identify the monodromy of a strictly linearly fibered 
arrangement.  Suppose $\A$ is strictly linearly fibered over $\B$, and 
write $|\B|=m$ and $|\A|=m+n$.  It then follows from the definition 
that a defining polynomial for $\A$ factors as
\begin{equation} \label{eq:defpoly}
Q(\A) = Q(\B)\cdot \phi(\bx,z)
\end{equation} 
where $Q(\B)=Q(\B)(\bx)$ is a defining polynomial for $\B$, and 
$\phi(\bx,z)$ is a product of $n$ linear functions:
\[
\phi(\bx,z) = (z-g_1(\bx))(z-g_2(\bx))\cdots (z-g_n(\bx)), 
\quad g_j(\bx)\ \text{linear.}
\]
Since $\phi(\bx,z)$ has distinct roots for any $\bx\in M(\B)$, the 
associated {\em root map}
\begin{equation} \label{eq:rootmap}
g:M(\B) \to \C^n,\qquad 
g(\bx)=\bigl(g_1(\bx),g_2(\bx),\dots,g_n(\bx)\bigr),
\end{equation}
takes values in the configuration space $F(\C,n)$.

\begin{thm}\label{thm:slfmono}
Let $\B$ be an arrangement of $m$ hyperplanes, and let $\A$ be an 
arrangement of $m+n$ hyperplanes which is strictly linearly fibered 
over $\B$.  Then the bundle $p:M(\A)\to M(\B)$ is equivalent to the 
pullback of the bundle of configuration spaces $p_{n+1}:F(\C,n+1)\to 
F(\C,n)$ along the map $g$.
\end{thm}
\begin{proof}
Denote points in $F(\C,n+1)$ by $(\by,z)$, where 
$\by=(y_1,\dots,y_n)\in F(\C,n)$ and $z\in \C$ satisfies $z\neq y_j$ 
for each $j$.  Similarly, denote points in $M(\A)$ by $(\bx,z)$, where 
$\bx\in M(\B)$ and $\phi(\bx,z)\neq 0$.  Then $p_{n+1}(\by,z)=\by$ and 
$p(\bx,z)=\bx$.

Let $E=\bigl\{\bigl(\bx,(\by,z)\bigr) \in M(\B) \times F(\C,n+1) \mid 
g(\bx)=\by\bigr\}$ be the total space of the pullback of 
$p_{m+1}:F(\C,n+1)\to F(\C,n)$ along $g$.  It is then readily checked 
that the map $M(\A) \to E$ defined by $(\bx,z) \mapsto 
\bigl(\bx,(g(\bx),z)\bigr)$ is an equivalence of bundles.
\end{proof}

Denote the fundamantal group of the complement of an arrangement $\A$ 
by $G(\A)$, or simply by $G$ if the underlying arrangement is clear.  
We record some immediate consequences of the above result.

\begin{cor} \label{cor:slfgroup} 
Let $\B$ be an arrangement of $m$ hyperplanes, and let $\A$ be an 
arrangement of $m+n$ hyperplanes which is strictly linearly fibered 
over $\B$.  Then,
\begin{enumerate}
\item the bundle $p:M(\A) \to M(\B)$ admits a section; 

\item the structure group of the bundle $p:M(\A)\to M(\B)$ is the pure 
braid group $P_n$;

\item the monodromy of the bundle $p:M(\A)\to M(\B)$ factors as 
$\eta=\a_n\circ \c$, where $\a_n:P_{n} \to \Aut(\F_{n})$ is the Artin 
representation and $\c=g_*:G(\B) \to P_{n}$ is the map on fundamental 
groups induced by $g$; and

\item if $\B$ is $K(\pi,1)$, then $\A$ is $K(\pi,1)$ and the group 
$G(\A)$ is isomorphic to the semidirect product $\F_n \rtimes_\eta 
G(\B)$.
\end{enumerate}
\end{cor}

\begin{rem} \label{rem:slfpres}
Let $g_i$ be the homotopy class of a meridional loop about the 
hyperplane $H_i$ of $\B$.  The classes $g_i$ generate the group 
$G(\B)$, see \cite{OT}.  Suppose $\A$ is strictly linearly fibered 
over $\B$, and identify the fundamental group of the fiber of 
$p:M(\A)\to M(\B)$ with the free group $\F_n=\langle 
t_1,\dots,t_n\rangle$.  If $\B$ is $K(\pi,1)$, then so is $\A$, and 
$1\to \F_n \to G(\A) \leftrightarrows G(\B) \to 1$ is split exact.  
Identify $\F_n$ and $G(\B)$ with their images in $G(\A)$ under the 
inclusion and splitting respectively.  Then the group $G(\A)$ has 
presentation
\[
G(\A)=\langle g_i\ (1\le i\le m), t_j\ (1\le j\le n) \mid 
g_i^{-1} t_{j}^{} g_{i}^{} = 
\eta_{}^{}(g_{i}^{})(t_{j}^{})\ (1\le j\le n, 1\le i\le m)\rangle.
\]
The automorphism $\eta(g_{i})$ is obtained by applying the Artin 
representation \eqref{eq:ArtinRep} to the pure braid $\c(g_{i})$.  A 
method for calculating the braids $\c(g_{i})$ is presented in Section 
\ref{subsec:monocalc}.
\end{rem}

Now let $\A$ be a fiber-type arrangement in $\C^\ll$ with exponents 
$\{d_1,\dots,d_\ll\}$.  Then there is a choice of coordinates 
$(x_{1},\dots,x_{\ll})$ on $\C^{\ll}$ so that a defining polynomial 
for $\A$ factors as $Q(\A)=\prod_{k=1}^{\ll} Q_{k}(x_{1},\dots,x_{k})$, 
where $\deg Q_{k}=d_{k}$ for each $k$, see \eqref{eq:defpoly}.  For 
each $j\le \ll$, the polynomial $\prod_{k=1}^{j} Q_{k}$ defines a 
fiber-type arrangement $\A_{j}$ in $\C^{j}$ with exponents 
$\{d_1,\dots,d_j\}$, and $\A_{j}$ is strictly linearly fibered over 
$\A_{j-1}$.

Identify the fundamental group of the fiber of the bundle 
$p_j:M(\A_j)\to M(\A_{j-1})$ with the free group $\F_{d_j}$ on $d_j$ 
generators for each $j$.  The action of the group $G(\A_{j-1})$ on 
$\F_{d_j}$ is the composition, $\eta_j:=\a_{d_j} \circ \c_j$, of the 
Artin representation $\a_{d_j}:P_{d_j} \to \Aut(\F_{d_j})$ and the 
homomorphism $\c_j:G(\A_{j-1}) \to P_{d_j}$ induced by the map $g_j$ 
of \eqref{eq:rootmap}.  Repeated application of 
Theorem~\ref{thm:slfmono} and Corollary~\ref{cor:slfgroup} yields

\begin{thm} \label{thm:ftmono}
The fundamental group $G(\A_\ll)$ of the complement of the fiber-type 
arrangement $\A_{\ll}$ admits the structure of an iterated semidirect 
product of free groups
\[
G(\A_\ll) = \F_{d_\ll} \rtimes_{\eta_{\ll}} \cdots \rtimes
\F_{d_2} \rtimes_{\eta_{2}}  \F_{d_1},
\]
with (split) pure braid extensions $\eta_j:G(\A_{j-1})\to P_{d_j} < 
\Aut(\F_{d_j})$.
\end{thm}

\begin{rem} \label{rem:ftpres}
For $1\le j\le \ll$, fix generators $x_{q,j}$, $1\le q \le d_j$, for 
the free group $\F_{d_j}$.  Then the group $G=G(\A_\ll)$ has 
presentation
\[
G=\langle x_{q,j}^{}\ (1\le q \le d_j,\ 1\le j \le \ll) \mid 
x_{p,i}^{-1} x_{q,j}^{} x_{p,i}^{} = 
\eta_j^{}(x_{p,i}^{})(x_{q,j}^{})\ (i<j)\rangle.
\]
\end{rem}

\subsection{Fiber-Type Arrangement Groups are Strongly Poly-Free} 
\label{subsec:SPF}
\  
\medskip

As noted in the Introduction, the structure of the fundamental group 
of the complement of a fiber-type arrangement exhibited in Theorem 
\ref{thm:ftmono} may be used to obtain a number of interesting and 
important consequences.  In this section we record another, showing 
that these groups are strongly poly-free.  We first recall the 
definition of this class of groups from \cite{AFR}.

\begin{defn} \label{def:SPF}
A discrete group $G$ is {\em strongly poly-free} if there exists a 
finite filtration of $G$ by subgroups, $1=G_0 < G_1 < \dots < G_\ll 
=G$, which satisfies the conditions:
\begin{enumerate}
\item[(1)] $G_k$ is normal in $G$ for each $k$;
\item[(2)] $G_{k+1}/G_k$ is a finitely generated free group; and
\item[(3)] for each $w\in G$ and each $k$, there is a compact surface 
$\Sigma$ 
with non-empty boundary
and a diffeomorphism $f:\Sigma\to\Sigma$ such that the 
induced homomorphism $f_*$ on $\pi_1(\Sigma)$ is equal to $\conj_w$ 
in $\Out(\pi_1(\Sigma))$, where $\conj_w$ is the action of $w$ on 
$G_{k+1}/G_k$ by conjugation and $\pi_1(\Sigma)$ is identified with 
$G_{k+1}/G_k$ via a suitable isomorphism.
\end{enumerate}
\end{defn}

\begin{thm}[Aravinda, Farrell, and Roushon \cite{AFR}, Theorem~2.1]
\label{thm:purebraidSPF}
For each $n$, the pure braid group $P_{n}$ is strongly poly-free.
\end{thm}

In light of this result, and those of Section~\ref{subsec:ftmono}, it 
is natural to speculate that the fundamental group of the complement 
of any fiber-type arrangement is strongly poly-free.  This is indeed 
the case.

\begin{thm} \label{thm:ftSPF}
Let $\A=\A_{\ll}$ be a fiber-type arrangement.  Then the fundamental 
group $G=G(\A_\ll)$ of the complement is strongly poly-free.
\end{thm}
\begin{proof}
 From Theorem \ref{thm:ftmono}, we have $G\cong \F_{d_\ll} 
\rtimes_{\eta_{\ll}} \cdots \rtimes_{\eta_{3}} \F_{d_2} 
\rtimes_{\eta_{2}} \F_{d_1}$.  For $1 \le k \le \ll$, let 
$G_{k}=\F_{d_\ll} \rtimes_{\eta_{\ll}} \cdots \rtimes_{\eta_{\ll-k+2}} 
\F_{d_{\ll-k+1}}$.  Then $G_{k}$ is normal in $G$ and 
$G_{k+1}/G_{k}=\F_{d_{\ll-k}}$ is free, so conditions (1) and (2) of 
the definition are satisfied.

We show that condition (3) holds by induction on the cohomological 
dimension of $G$, which we may assume without loss of generality is 
equal to $\ll$.  In the case $\ll=1$, $G=\F_{d}$ is a (single) 
finitely generated free group acting on itself by conjugation, and 
condition (3) clearly holds.  In general, from Theorem 
\ref{thm:slfmono} we have a commuting diagram
\[
\begin{CD}
1 @>>> \F_{d_{\ll}} @>>> G(\A_{\ll}) @>>> G(\A_{\ll-1}) @>>>1\\
@.     @VV{\id}V   @VV{\tilde\c_{\ll}}V  @VV{\c_{\ll}}V\\
1 @>>> \F_{d_{\ll}} @>>> P_{d_{\ll}+1} @>>> P_{d_{\ll}} @>>>1 
\end{CD}
\]
where $\c_\ll$ is induced by the map $g_\ll:M(\A_{\ll-1})\to 
F(\C,d_\ll)$ from \eqref{eq:rootmap}, and $\tilde\c_\ll$ is induced 
by $\tilde g_\ll:M(\A_{\ll})\to F(\C,d_\ll+1)$ defined by $\tilde 
g_\ll(\bx,z)=\bigl(g_\ll(\bx),z\bigr)$.  From 
Corollary~\ref{cor:slfgroup}, we have $G(\A_\ll)=\F_{d_{\ll}} 
\rtimes_{\eta_\ll} G(\A_{\ll-1})$, where 
$\eta_\ll = \a_{d_\ll} \circ \c_\ll$ and $\a_{d_\ll}:P_{d_{\ll}} \to 
\Aut(\F_{d_{\ll}})$ is the Artin representation.  Using the 
semidirect product structure, every element $w \in G(A_\ll)$ may be 
expressed as $w=uv$, where $u\in\F_{d_\ll}$ and $v \in G(\A_{\ll-1})$.

For $w \in G(\A_\ll)$, consider the conjugation action of $w$ on 
$G_{k+1}/G_k$.  In the case $k=0$, for $a\in G_1=\F_{d_\ll}$ and 
$w=uv$ as above, we have
\[
w^{-1} a w = v^{-1} u^{-1} a u v = \eta_\ll(v)(u^{-1} a u) = 
\c_\ll(v)^{-1} u^{-1} a u \c_\ll(v).
\]
Thus in this instance, conjugation by $w$ coincides with conjugation 
by the pure braid $u\cdot \c_\ll(v) =\tilde\c_\ll(uv)=\tilde\c_\ll(w) 
\in P_{d_\ll+1}$.  So for $k=0$, condition (3) holds by the result of 
Aravinda, Farrell, and Roushon stated in Theorem 
\ref{thm:purebraidSPF} above.

For the case $k>0$, let $a \in G(\A_{\ll-1})$ and consider $w^{-1}aw$.  
In this instance, we have
\[
w^{-1} a w = v^{-1} u^{-1} a u v = v^{-1}(u^{-1} a u a^{-1})v \cdot 
v^{-1} a v.
\]
Now $v^{-1}(u^{-1} a u a^{-1})v$ is in $G_1$, since $u \in G_1$ and 
$G_1$ is normal in $G$.  Consequently, the conjugation action of 
$w=uv$ on $G_{k+1}/G_k$ coincides with that of $v\in G(\A_{\ll-1})$ 
for $k>0$.  So condition (3) holds by induction in this case.
\end{proof}

In \cite[Theorem 1.3]{AFR}, it is shown that the Whitehead group of 
the direct product of a strongly poly-free group $G$ and a free 
abelian group is trivial, $Wh(G\times\Z^m)=0$ for every $m\ge 0$.  
This result and Theorem \ref{thm:ftSPF} above yield the following, 
which resolves positively the conjecture of Aravinda, Farrell, and 
Roushon stated in \cite[Section 2]{AFR}.

\begin{cor} \label{cor:ftWh}
Let $G$ be the fundamental group of the complement of a fiber-type 
arrangement $\A$.  Then the Whitehead group of $G\times\Z^{m}$ is 
trivial for every $m\ge 0$.
\end{cor}

\begin{rem} 
For $G=G(\A)$ as above, note that the group $G\times\Z^{m}$ may itself 
be realized as the fundamental group of the complement of the 
fiber-type arrangement $\A\times\B$, where $\B$ is the arrangement of 
coordinate hyperplanes in $\C^m$.
\end{rem}

\subsection{Calculating the Monodromy} \label{subsec:monocalc}
\ 
\medskip

In this section, we present a method for calculating the monodromy of 
the bundle $p:M(\A)\to M(\B)$ for an arrangement $\A$ of $m+n$ 
hyperplanes, strictly linearly fibered over the arrangement $\B$ of 
$m$ hyperplanes.  This technique may be applied repeatedly to 
determine the iterated semidirect product structure of the group of a 
fiber-type arrangement, as we illustrate in the next section.  Since 
the monodromy factors as $\eta= \a_n \circ \c:G(\B) \to P_n \to 
\Aut(\F_n)$, where $\a_n:P_n \to \Aut(\F_n)$ is the Artin 
representation, we focus on the determination of the homomorphism 
$\c:G(\B)\to P_n$ induced by the map $g:M(\B)\to F(\C,m)$ from 
\eqref{eq:rootmap}.  As this is a natural generalization of the method 
for finding the braid monodromy of a complex line arrangement 
developed in \cite{CSbm}, we call $\c:G(\B) \to P_n$ the {\em braid 
monodromy} of the bundle $p:M(\A)\to M(\B)$.

Write $\B=\{H_1,\dots,H_m\}$, and let $L$ be a complex line in 
$\C^\ll$ that is transverse to $\B$.  Denote the coordinate on $L$ by 
$x$, and the point $L\cap H_j$ by $q_j$.  Then $C:=L \setminus 
\{q_1,\dots,q_m\}$ is the complement of $m=|\B|$ points in $L=\C$, and 
$\pi_{1}(C)=\F_{m}$.  Let $i:C\hookrightarrow M(\B)$ denote the 
natural inclusion, and let $\widehat g=g\circ i:C\to F(\C,n)$ denote 
the restriction to $C$ of the map $g:M(\B)\to F(\C,n)$.  Passing to 
fundamantal groups, we have $\widehat{g}_*=\widehat\c=\c \circ 
i_{*}:\F_{m} \to G(\B)\to P_{n}$.  Since $i_{*}$ is surjective, it 
suffices to determine the homomorphism $\widehat\c:\F_{m}\to P_{n}$.

The pullback of the configuration space bundle $F(\C,n+1)\to F(\C,n)$ 
along $\widehat{g}$ is equivalent to the restriction, 
$\widehat{p}:Y\to C$, of the bundle $p:M(\A)\to M(\B)$, where
\[
Y = \{(x,z) \in C\times \C \mid \phi(x,z) \neq 0\}=
\{(x,z)\in\C^2 \mid \prod_{j=1}^n(x-q_j) \cdot \phi(x,z) \neq 0\},
\]
and $\phi(x,z)$ is the restriction of $\phi(\bx,z)$ to $L$.  The 
polynomial $\phi(x,z)$ defines an arrangement $\H$ of $n$ lines in 
$\C^2$.  The multiple points of $\H$ necessarily lie on the lines 
$x=q_j$.  Note that more than one such multiple point may lie on a 
given such line, and that there may be lines $x=q_j$ upon which no 
such multiple points lie.  The present construction generalizes that 
of \cite{CSbm} in these senses.

Order the $m$ distinct points $q_{j}$ in $L=\C$ by decreasing real 
part, breaking ties by imaginary part.  If $i<j$, then $\Real(q_{i}) > 
\Real(q_{j})$, or $\Real(q_{i}) = \Real(q_{j})$ and $\Imag(q_{i}) < 
\Imag(q_{j})$.  Fix a basepoint $q_0$ in $C$ with 
$\Real(q_{0})>\Real(q_{1})$, and let $\xi=\xi(t)$ be a path in $\C$, 
emanating from $q_0$ and passing through the ordered points $q_j$.  In 
a small disk $D_\epsilon(q_j)$ about $q_j$, take $\xi$ to be a 
horizontal line segment, which passes through $q_j$ from right to left 
as $t$ increases, and let $q'_j=q_j-\epsilon$ and 
$q_j''=q_j+\epsilon$.  Choose $\epsilon>0$ small, so that $q_i \notin 
D_\epsilon(q_j)$ for $i \neq j$.  Let $\xi_{j,j+1}$ be the portion of 
$\xi$ from $q_j''$ to $q_{j+1}'$, let $\xi_j'$ be the portion of the 
boundary of the disk $D_\epsilon(q_j)$ from $q_j'$ to $q_j''$, and 
$\xi_j''$ the portion of ${\partial}D_\epsilon(q_j)$ from $q_j''$ to 
$q_j'$ (both oriented counterclockwise).

Let $u_{j}$ denote the homotopy class of the loop in $C$ based at 
$q_{0}$ which traverses the paths $\xi_{i,i+1}$ and $\xi_i'$ for $i<j$ 
in the natural order, passes around $q_j$ along $\xi_j'$ and 
$\xi_j''$, and returns to $q_0$ along the $\xi_{i,i+1}$ and $\xi_i'$ 
with $i<j$.  Using these meridians, identify $\pi_1(C,q_0)$ with 
$\F_m=\langle u_{1},\dots,u_{m}\rangle$.  The monodromy of the bundle 
$\widehat{p}:Y \to C$ is determined by the (pure) braids 
$\widehat\c(u_{j})$.  Since the images $v_j=i_*(u_j)$ generate 
$G(\B)=\pi_1(M(\B))$, these braids also determine the monodromy of 
$p:M(\A) \to M(\B)$.

The braids $\widehat\c(u_{j})=\c(v_{j})$ may be calculated from the 
{\em braided wiring diagram} $\W = \{(x,z)\in \xi\times\C \mid 
\phi(x,z)=0\}$ associated to the path $\xi$, cf.~\cite[Section 
5]{CSbm}.  For $q \neq q_{j}$, $1\le j \le m$, the set $\W\cap\{x=q\}$ 
consists of $n$ distinct points, the intersections of the lines of 
$\H$ with $\{x=q\}$.  Order the lines of $\H$ by increasing real part 
of the $n$ points of $\W\cap\{x=q_{0}\}$, breaking ties as above.  Let 
$[n]=\{1,\dots,n\}$.  The diagram $\W$ may be (abstractly) specified 
by a sequence of partitions of $[n]$ and braids,
\[
\W=\W_m=
\{I(1),\b_{1,2},I(2),\b_{2,3},\dots,\b_{m-1,m},I(m),\b_{m,m+1}\}.
\]
The braids $\b_{i,i+1}$ are elements of the full braid group $B_n$, 
obtained by tracing the components of $\W$ over the path 
$\xi_{i,i+1}$, see \cite[Section 4.4]{CSbm}.  The partitions
\[
I(j)=\bigl(I_{1}(j) \mid I_{2}(j) \mid \dots \mid I_{r}(j)\bigr)=
\bigl(1,\dots,j_{1}\mid j_{1}+1,\dots,j_{2} \mid \cdots \cdots \mid
j_{r-1}+1,\dots,n\bigr)
\]
record the ordering at $x=q_j'$ of the lines of $\H$ which meet at 
$x=q_{j}$.  If no lines of $\H$ meet at $x=q_j$, then $I(j)=(1\mid 
2\mid\dots\mid n)$ consists of $n$ singletons.

To each block $I_{k}=\{i,i+1,\dots,i+s\}$ of such a partition $I$, we 
associate a permutation braid $\mu_{I_{k}}$, a half twist on $I_{k}$, 
given in terms of the standard generators of $B_{n}$ by
\[
\mu_{I_{k}}=(\s_{i}\cdots\s_{i+s-1})(\s_{i}\cdots\s_{i+s-2})\cdots
(\s_{i}\cdots\s_{i+1})(\s_{i}).
\]
If $|I_k|=1$, set $\mu_{I_k}=1$.  Note that $\mu_{I_{k}}$ and 
$\mu_{I_{k'}}$ commute for $k\neq k'$, and that the product 
$\Upsilon_{I}=\mu_{I_{1}} \cdots \mu_{I_{r}}$ records the braiding of 
the components of $\W$ over each of the paths $\xi_{j}'$ and 
$\xi_{j}''$.  The local (braid) monodromy around the point $q_{j}$ is 
given by $\Upsilon_{I}^{2}$, the product of the full twists 
$\mu_{I_{k}}^{2}$.  The braid monodromy $\widehat{\c}:\F_{m}\to P_{n}$ 
is then given by
\begin{equation} \label{eq:bmfull}
\widehat\c(u_{j}) = \b_{j}^{-1} \Upsilon_{I(j)}^{2} \b_{j}^{}, 
\end{equation}
where the conjugating braids $\b_{j}$ satisfy $\b_{1}=1$ and 
$\b_{j+1}=\b_{j,j+1} \cdot \Upsilon_{I(j)} \cdot \b_{j}$ for $j\ge 
1$.  

We express the braid monodromy solely in terms of pure braids, 
cf.~\cite[Section~5.3]{CSbm}.  Recall the original ordering of the 
lines of $\H$ at the basepoint $x=q_{0}$.  Let 
$V(j)=\bigl(V_{1}(j)\mid V_{2}(j)\mid \dots \mid V_{r}(j)\bigr)$ be 
the partition of $[n]$ recording the indices of these lines meeting at 
$x=q_{j}$ in terms of this ordering.  To a block 
$V_{k}=\{k_{1},\dots,k_{s}\}$ of such a partition (with 
$k_{i}<k_{i+1}$), associate the full twist on $V_{k}$, given in terms 
of the standard generators of $P_{n}$ by
\[
A_{V_k} = (A_{k_1,k_2})\cdot (A_{k_1,k_3}A_{k_2,k_3})\cdots\cdots
(A_{k_1,k_s}\cdots A_{k_{s-1},k_s}).
\]
Geometrically, these braids are obtained by gathering the strands 
indexed by $V_k$ together behind the remaining strands, performing a 
full twist on the $V_k$ strands, and then returning these strands to 
their original positions.  If $s=1$, set $A_{V_{k}}=1$.  Expressing 
\eqref{eq:bmfull} in terms of pure braids yields

\begin{thm} \label{thm:bmpure}
The braid monodromy $\widehat{\c}:\F_m\to P_n$ of the bundle 
$\widehat{p}:Y\to C$, and hence that of 
the strictly linearly fibered bundle 
$p:M(\A) \to M(\B)$, is given by
\begin{equation} \label{eq:bmpure}
\c\circ i_*(u_j)=\widehat\c(u_{j}) = 
\prod_{k=1}^r A_{V_k(j)}^{\zeta_k} = \prod_{k=1}^r
\zeta_k^{-1}\cdot A_{V_k(j)}\cdot{\zeta_k},
\end{equation}
where $\zeta_k\in P_{n}$ is determined by the subdiagram 
$\W_{j-1}$ of $\W$ and the block $V_k(j)$.
\end{thm}

\begin{rem} \label{rem:mono&comb}
If $\A$ is strictly linearly fibered over $\B$, the associated braid 
monodromy determines the rank two elements of the intersection poset 
$L(\A)$ in the following way.  Order the hyperplanes of 
$\A\setminus\B=\H=\{\H_1,\dots,\H_n\}$ as indicated above.  Let $g\in 
G(\B)$ be a meridian about a hyperplane $H$ of $\B$, with associated 
monodromy generator $\c(g) = \prod_{k=1}^r A_{V_k}^{\zeta_k}$.  Then 
the blocks $V_k$ record those hyperplanes of $\A\setminus\B$ which 
meet $H$ in codimension two.  In other words, for each $k$, we have $H 
\cap \bigl(\bigcap_{j\in V_k} \H_j \bigr) \in L_2(\A)$.  Thus, 
$L_2(\A)$ consists of elements of this form, together with elements of 
$L_2(\B)$.
\end{rem}

In addition to recording combinatorial information, the expression 
\eqref{eq:bmpure} of the braid monodromy also sheds light on the 
presentations of strictly linearly fibered and fiber-type arrangement 
groups noted in Remarks \ref{rem:slfpres} and \ref{rem:ftpres}.  To 
this end, we briefly describe the behavior of pure braids of the form 
\eqref{eq:bmpure} under the Artin representation.  Note that the 
factors of these braids commute, 
$[A_{V_k(j)}^{\zeta_k},A_{V_l(j)}^{\zeta_l}]=1$.  For $t \in \F_n$, 
let $\widetilde{t}$ denote some conjugate of $t$.  Let $V=(V_{1}\mid 
V_{2}\mid \dots \mid V_{r})$ be a partition of $[n]$.  If 
$V_k=\{k_1,\dots,k_s\}$ is a block of $V$, set 
$\widetilde{t}_{V_k}=\widetilde{t}_{k_1}\cdots \widetilde{t}_{k_s}$.

\begin{prop} \label{prop:bmartin}
Let $\c=\prod_{k=1}^r A_{V_k}^{\zeta_k}$ be a braid with commuting 
factors associated to the partition $V$ of $[n]$, and let 
$\a_n:P_n\to\Aut(\F_n)$ be the Artin representation.  Then
\begin{equation*} \label{eq:bmartin}
\a_n(\c)(t_j) = \widetilde{t}_{V_k}^{} \cdot \widetilde{t}_{j}^{} 
\cdot 
\widetilde{t}_{V_k}^{-1},
\end{equation*}
where $V_k$ is the unique block of $V$ containing $j$.
\end{prop}
\begin{proof}
Write $\c=B\cdot A_{V_k}^{\zeta_k}$, where $j\in V_k$ and 
$B=\prod_{l\neq k} A_{V_l}^{\zeta_l}$, and identify a pure braid with 
its image under the Artin representation.  Then 
$\a_n(\c)(t_j)=\c(t_j)=A_{V_k}^{\zeta_k} \circ B(t_j)$.  The action of 
$P_n<\Aut(\F_n)$ is by conjugation, so $\c(t_j)=\zeta_k\circ 
A_{V_k}(w\cdot t_j\cdot w^{-1})$, where $w\cdot t_j\cdot 
w^{-1}=\zeta_k^{-1}\circ B(t_j)$.  A calculation with the Artin 
representation reveals that $A_{V_k}(t_j)= t_{V_k}^{} \cdot t_{j}^{} 
\cdot t_{V_k}^{-1}$ for $j\in V_k$.  So we have $\c(t_j)=v \cdot 
\z_k(t_{V_k}^{} \cdot t_{j}^{} \cdot t_{V_k}^{-1}) \cdot v^{-1}$, 
where $v=\z_k\circ A_{V_k}(w)$, and the result follows.
\end{proof}

\subsection{The Coxeter Arrangement of type $\BB$} 
\label{subsec:typeB}
\ 
\medskip

Let $\B_{n}$ denote the Coxeter arrangement of type $\BB$ in $\C^{n}$, 
with defining polynomial
\[
Q(\B_{n})=x_{1}\cdots x_{n} \prod_{i<j}(x_{j}^{2}-x_{i}^{2}),
\]
and complement $M(\B_{n})$.  As shown by Brieskorn \cite{Br}, this 
arrangement is fiber-type, and the fundamental group of $M(\B_{n})$ is 
the pure braid group of type $\BB$, $G(\B_n)=PB_{n}$.  We illustrate 
the method described in the previous section by determining the 
iterated semidirect product structure of this generalized pure braid 
group.

Denote the $n^{2}$ hyperplanes of $\B_{n}$ by $H_{i}=\ker(x_{i})$ and 
$H_{i,j}^{\pm}=\ker(x_{j}\pm x_{i})$.  The line 
$L=\{(x,2x+b_{2},\dots,nx+b_{n})\}$, where $b_{k}=(2k+1)!$, is 
transverse to $\B_{n}$.  Write $L\cap H_{i}=h_{i}$ and $L\cap 
H_{i,j}^{\pm}=h_{i,j}^{\pm}$.  Notice that these points are real, and 
check that
\[
h_{n-1,n}^{-}<\dots < h_{1,n}^{-} < h_{n} < h_{1,n}^{+}<\dots < 
h_{n-1,n}^{+}<\cdots \cdots  <
h_{1,2}^{-} < h_{2} < h_{1,2}^{+} < h_{1}.
\]
With this notation, we have $C=L \setminus 
\{h_i\}\cup\{h_{i,j}^\pm\}$.

\begin{figure}[h]
\epsfysize=1.5 truein \epsfbox{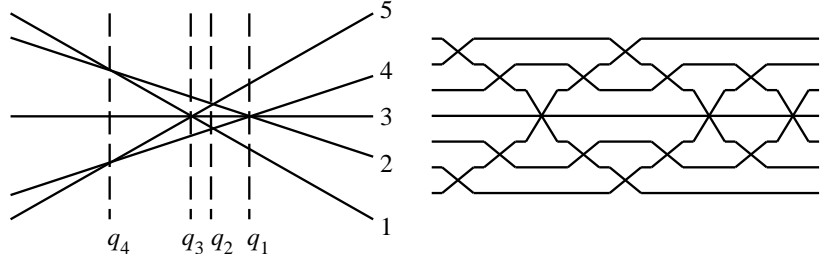}
\caption{Type $\BB$ wiring diagrams for $n=2$ (left) and $n=3$ 
(right)}
\label{fig:BW}
\end{figure}

The arrangement $\B_{n+1}$ is strictly linearly fibered over $\B_n$, 
and has defining polynomial $Q(\B_{n+1})=Q(\B_n) \cdot \phi_n(\bx,z)$, 
where
\[
\phi_n(\bx,z)=(z+x_n)\cdots(z+x_1)\cdot z \cdot (z-x_1) \cdots 
(z-x_n).
\]
Let $g_n:M(\B_n)\to F(\C,2n+1)$ be the associated root map 
(cf.~\eqref{eq:rootmap}), inducing $\c_n:PB_n \to P_{2n+1}$ on 
fundamental groups, and let $\widehat{g}_n$ denote the restriction of 
$g_n$ to $C$.  The restriction, $\phi_n(x,z)$, of $\phi_n(\bx,z)$ to 
$L$ defines an arrangement $\H^n$ of $2n+1$ lines in $\C^2$.  These 
lines have real defining equations, so let $q_0 \in C$ be a real 
basepoint with $h_1<q_0$, and let $\xi=[d,q_0]$ be a line segment 
along the real axis in $L$ with $d < h_{n-1,n}^-$.  The resulting 
wiring diagram $\B\W^n$ is unbraided: the braids $\b_{i,i+1}=1$ are 
all trivial.  The arrangement $\H^2$ and diagram $\B\W^3$ are depicted 
in Figure~\ref{fig:BW}.

\begin{exm} \label{exm:B2}
We explicitly carry out the monodromy calculation in the case $n=2$.  
Write $q_1=h_1$, $q_2=h_{1,2}^-$, $q_3=h_2$, and $q_4=h_{1,2}^+$.  
Refering to Figure \ref{fig:BW}, we see that the partitions $I(j)$ 
and associated braids $\Upsilon_{I(j)}$ are
\[
\begin{alignat*}{5}
I(1)&=I(3)&=(1\mid 2,3,4 \mid 5),&\qquad\qquad
\Upsilon_{I(1)}\,&=\Upsilon_{I(3)}&=1\cdot \s_2\s_3\s_2\cdot 
1,\qquad\\
I(2)&=I(4)&=(1,2\mid 3 \mid 4,5),&\qquad\qquad
\Upsilon_{I(2)}\,&=\Upsilon_{I(4)}&=\s_1 \cdot 1 \cdot \s_4.\qquad
\end{alignat*}
\]
Thus the braids $\b_j$ appearing in \eqref{eq:bmfull} are
$\b_1=1$, $\b_2=\s_2\s_3\s_2$, $\b_3=\s_1\s_4\s_2\s_3\s_2$, and
$\b_4=\s_2\s_3\s_2\s_1\s_4\s_2\s_3\s_2$.

Write $\pi_1(C)=\F_4=\langle u_1,u_2,u_3,u_4\rangle$ as before.  
Combing the braid as in \eqref{eq:bmpure}, the braid 
monodromy $\widehat\c_2:\F_4 \to P_5$ is given by
\[
\begin{alignat*}{3}
\widehat\c_2(u_1)&=A_{2,3,4}^{}, \qquad\qquad
&\widehat\c_2(u_2)&=A_{1,2}^{}A_{1,3}^{}A_{1,4}^{}
A_{1,3}^{-1}A_{1,2}^{-1} \cdot A_{2,5}^{},\qquad\\
\widehat\c_2(u_3)&=A_{1,2}^{}A_{1,3,5}^{}A_{1,2}^{-1},\qquad\qquad
&\widehat\c_2(u_4)&=A_{1,2}^{} \cdot A_{4,5}^{}.\qquad
\end{alignat*}
\]
These braids act on $\F_5=\langle t_1,\dots,t_5\rangle$ via the Artin 
representation $\a_5:P_5\to\Aut(\F_5)$.  Write 
$\widehat\eta_2=\a_5\circ\widehat\c_2$.  
A calculation using \eqref{eq:ArtinRep} yields
$\widehat\eta_2(u_i)(t_j)=w_{i,j}^{} t_j w_{i,j}^{-1}$, where
\[
\begin{alignat*}{2}
w_{1,j}&=\begin{cases}
1,&\text{if $j=1,5$,}\\
t_2t_3t_4&\text{if $j=2,3,4$,}
\end{cases} \qquad
&w_{2,j}&=\begin{cases}
t_1^{}t_2^{}t_3^{}t_4^{}t_3^{-1}t_2^{-1}&\text{if $j=1$,}\\
t_2t_5&\text{if $j=2,5$,}\\
{[}t_2,t_5{]}&\text{if $j=3$,}\\
{[}t_2,t_5{]}t_3^{-1}t_2^{-1}t_1^{}t_2^{}t_3^{}&\text{if $j=4$,}
\end{cases} \qquad\\
w_{3,j}&=\begin{cases}
t_1^{}t_2^{}t_3^{}t_5^{}t_2^{-1}&\text{if $j=1$,}\\
1,&\text{if $j=2$,}\\
t_2^{-1}t_1^{}t_2^{}t_3^{}t_5^{}&\text{if $j=3,5$,}\\
{[}t_2^{-1}t_1^{}t_2^{}t_3^{},t_5^{}{]}&\text{if $j=4$,}
\end{cases} \qquad
&w_{4,j}&=\begin{cases}
t_1t_2&\text{if $j=1,2$,}\\
1&\text{if $j=3$,}\\
t_4t_5&\text{if $j=4,5$.}
\end{cases} \qquad
\end{alignat*}
\]

The braid monodromy $\widehat\c_2:\F_4\to P_5$ descends to 
$\c_2:PB_2\to P_5$.  With the above Artin representation calculation, 
this realizes the group $PB_3$ as a semidirect product, $PB_3 = \F_5 
\rtimes_{\eta_2} PB_2$.  We momentarily defer further discussion of 
this realization.
\end{exm}

In general, the group $PB_n$ is generated by the images of the 
generators $u_j$ of $\pi_1(C)$ under the map induced by inclusion $i:C 
\hookrightarrow M(\B_n)$.  These images, $c_j$, $a_{i,j}$, and 
$b_{i,j}$, are homotopy classes of meridional loops about the 
hyperplanes $H_j$, $H_{i,j}^-$, and $H_{i,j}^+$ respectively.  So for 
instance, $c_1=i_*(u_1)$, $b_{1,2}=i_*(u_2)$, $c_2=i_*(u_3)$, and 
$a_{1,2}=i_*(u_4)$.  The structure of the wiring diagram $\B\W^n$ is 
analogous to those exhibited in the examples shown in Figure 
\ref{fig:BW}.  Analysis of this structure and the ensuing braid 
monodromy calculation are left to the reader.  Define $U_{r,s}\in 
P_{2n+1}$ by $U_{r,s}=A_{r,r+1}A_{r,r+2}\cdots A_{r,s}$.  The result 
is

\begin{prop} \label{prop:Bbmono}
The braid monodromy $\c:PB_n \to P_{2n+1}$ is given by
\[
\begin{align*}
\c(c_j)&=U_{n-j+1,n}^{} A_{n-j+1,n+1,n+j+1}^{} U_{n-j+1,n}^{-1}\\
\c(a_{i,j})&=A_{n+i+1,n+j+1}^{} \cdot U_{n-j+1,n-i} 
A_{n-j+1,n-i+1}U_{n-j+1,n-i}^{-1},\\
\c(b_{i,j})&=A_{n-i+1,n+j+1}^{} \cdot U_{n-j+1,n+i} 
A_{n-j+1,n+i+1}U_{n-j+1,n+i}^{-1}.
\end{align*}
\]
\end{prop}

As in the case $n=2$ above, applying the Artin representation yields 
the semidirect product structure of $PB_{n+1}=\F_{2n+1} 
\rtimes_{\eta_n} PB_n$.  The correspondence between the generators 
$t_i$ of $\F_{2n+1}$ and the meridional generators of $PB_{n+1}$ is
\[
t_i=\begin{cases}
b_{n-i+1,n+1}&\text{if $1\le i \le n$,}\\
c_{n+1}&\text{if $i=n+1$,}\\
a_{i-n-1,n+1}&\text{if $n+2\le i \le 2n+1$.}
\end{cases}
\]
Carrying out the aforementioned Artin representation calculations for 
$k=1,\dots,n-1$ yields a presentation of the group $PB_{n}$ which 
exhibits the iterated semidirect product structure.  We suppress these 
lengthy calculations, and state the result below.  For $i<j$, let 
$\bar{a}_{i,j}=a_{1,j}a_{2,j}\cdots a_{i-1,j}$ and 
$\bar{b}_{i,j}=b_{i-1,j}\cdots b_{2,j}b_{1,j}$, and write 
$\hat{b}_{i,j}=b_{i,j}^{\bar{b}_{i,j}}$.

\begin{thm} \label{thm:PBn}
The Brieskorn generalized pure braid group admits the 
structure of an iterated semidirect product of free groups, 
$PB_n=\F_{2n-1} \rtimes \dots \rtimes \F_3 \rtimes \F_1$, 
where $\F_{2j-1}=\langle c_j,a_{i,j}, b_{i,j}\ (1\le i < j)\rangle$.  
For $j<l$, the action of $\F_{2j-1}$ on $\F_{2l-1}$ is given by
\[
\begin{alignat*}{3}
a_{i,j}^{-1}a_{k,l}^{} a_{i,j}^{}&=p_k a_{k,l} p_k^{-1},&\qquad
a_{i,j}^{-1}b_{k,l}^{} a_{i,j}^{}&=q_k b_{k,l} q_k^{-1},&\qquad
a_{i,j}^{-1}c_l^{} a_{i,j}^{}&=c_l,\\ 
b_{i,j}^{-1}a_{k,l}^{} b_{i,j}^{}&=u_k a_{k,l} u_k^{-1},&\qquad
b_{i,j}^{-1}b_{k,l}^{} b_{i,j}^{}&=v_k b_{k,l} v_k^{-1},&\qquad
b_{i,j}^{-1}c_l^{}b_{i,j}^{}&=w_l c_l w_l^{-1},\\ 
c_j^{-1}a_{k,l}^{} c_j^{}\,\ &=x_k a_{k,l} x_k^{-1},&\qquad
c_j^{-1}b_{k,l}^{} c_{s}^{}\,\ &=y_k b_{k,l} y_k^{-1},&\qquad
c_j^{-1}c_l^{} c_j^{}\,\ &=z c_l z^{-1},
\end{alignat*}
\]
where $w_l=[b_{i,l},a_{j,l}]$, 
$y_{j}=b_{j,l}\bar{b}_{j,l}c_{j}a_{j,l}\bar{b}_{j,l}^{-1}$, $y_{k}=1$
for $k\neq j$, 
$z^{}=\hat{b}_{j,l}c_la_{j,l}$ 
and
\[
\begin{align*}
p_k &= \begin{cases}
a_{i,l}a_{j,l}&\text{$k=i,j$,}\\
{[}a_{i,l},a_{j,l}{]}&\text{$i<k<j$,}\\
1&\text{else,}
\end{cases} \quad 
q_k = \begin{cases} 
\hat{b}_{j,l}^{\bar{b}_{r+1,j}^{-1}} 
&\text{$k=i$,}\\
b_{j,l}\hat{b}_{i,l}^{\bar{b}_{j,l}^{-1}}
&\text{$k=j$,}\\
1&\text{else,}
\end{cases} 
\quad 
x_k = \begin{cases} 
{[}\hat{b}_{j,l} c_l,a_{j,l}{]}
&\text{$k<j$,}\\
\hat{b}_{j,l} c_l
&\text{$k=j$,}\\
1&\text{else,}\\
\end{cases} \\
u_k &= \begin{cases} 
{[}b_{i,l},a_{j,l}{]}&\text{$k<i$ or $i<k<j$,}\\
{[}b_{i,l},a_{j,l}{]} \,
\hat{b}_{j,l}^{c_l \bar{a}_{i,l}}&\text{if $k=i$,}\\
b_{i,l}&\text{$k=j$,}\\
1&\text{else,}
\end{cases}
\quad
v_k = \begin{cases} 
{[}b_{i,l},a_{j,l}{]}&\text{$k<i$,}\\
b_{i,l} a_{j,l}&\text{$k=i$,}\\
b_{j,l}[\bar{b}_{j,l} c_l \bar{a}_{i,l}, a_{i,l}] a_{i,l}
&\text{$k=j$,}\\
1&\text{else.}
\end{cases}
\end{align*}
\]
\end{thm}

\begin{rem} In \cite[Section 3.8]{Le}, Leibman obtains a presentation 
of the group $PB_n$ by different means.  We have verified that this 
presentation and that of Theorem~\ref{thm:PBn} are equivalent for 
small $n$.  We have not been able to carry out this verification in 
general, as there appears to be a typographical error in \cite{Le} 
affecting the general case.
\end{rem}

\section{Orbit Configuration Spaces} \label{sec:orbit}

\subsection{Orbit Configuration Space Bundles} 
\label{subsec:OrbitBundle}
\ 
\medskip

Let $M$ be a manifold without boundary, and let $\Gamma$ be a finite 
group which acts freely on $M$.  The {\em orbit configuration space} 
consists of all ordered $n$-tuples of points in $M$ which lie in 
distinct orbits:
\[
F_\Gamma(M,n)=\{(x_1,\dots,x_n)\in M^n \mid 
\Gamma\cdot x_i\cap \Gamma\cdot x_j= \emptyset\ \text{if}\ i\neq j\}.
\]

Let $Q_n^\Gamma$ denote the union of $n$ distinct orbits, $\Gamma\cdot 
x_1,\dots,\Gamma\cdot x_n$, in $M$.  In \cite{Xi}, Xicot\'encatl 
proves the following theorem, a natural generalization to orbit 
configuration spaces of the Fadell-Neuwirth theorem stated in the 
Introduction.

\begin{thm}[Xicot\'encatl \cite{Xi}, Theorem~2.2.2]
For $\ll\le n$, the projection onto the first $\ll$ coordinates, 
$p_\Gamma:F_\Gamma(M,n)\to F_\Gamma(M,\ll)$, is a locally trivial 
bundle, with fiber $F_\Gamma(M\setminus Q_\ll^\Gamma,n-\ll)$.
\end{thm}

The proof given in \cite{Xi} is a natural adaptation of that of 
\cite{FN} for classical configuration spaces.  For the special case 
$n-\ll=1$, we give here a different proof, similar to that of Theorem 
\ref{thm:slfmono}, which sheds light on the structure of these 
bundles.

Suppose that the order of the finite group $\Gamma$ is $r$, and define 
a map from the orbit configuration space to the classical 
configuration space by sending a $n$-tuple of points in $M$ to their 
orbits.  Explicitly, define $f:F_\Gamma(M,k)\to F(M,rn)$ by 
$f(x_1,\dots,x_n)=(\Gamma\cdot x_1,\dots,\Gamma\cdot x_n)$.

\begin{thm} \label{thm:OrbitBundleThm}
The orbit configuration space bundle $p_\Gamma:F_\Gamma(M,n+1)\to 
F_\Gamma(M,n)$ is equivalent to the pullback of the bundle 
$p_{rn+1}:F(M,rn+1)\to F(M,rn)$ of classical configuration spaces 
along the map $f$.
\end{thm}
\begin{proof} Denote points in $F(M,rn+1)$ by $(\by,z)$, where 
$\by=(y_1,\dots,y_{rn})\in F(M,rn)$ and $z\in M$ satisfies $z\neq y_j$ 
for each $j$.  Similarly, denote points in $F_\Gamma(M,n+1)$ by 
$(\bx,z)$, where $\bm=(x_1,\dots,x_n)\in F_\Gamma(M,n)$ and $z\in M$ 
satisfies $\Gamma\cdot z \cap \Gamma\cdot x_j=\emptyset$ for each $j$.  
Then $p_{rn+1}(\by,z)=\by$ and $p_\Gamma(\bx,z)=\bx$.

Let $E=\{\bigl(\bx,(\by,z)\bigr)\in F_\Gamma(M,n)\times F(M,rn+1) \mid 
f(\bx)=\by\}$ be the total space of the pullback of 
$p_{rn+1}:F(M,rn+1)\to F(M,rn)$ along $f$.  It is then readily checked 
that the map $F_\Gamma(M,n+1)\to E$ defined by 
$(\bx,w)\mapsto\bigl(\bx,(f(\bx),w) \bigr)$ is an equivalence of 
bundles.
\end{proof}

We are mainly interested in the case where the finite cyclic group 
$\Gamma=\Z/r\Z$ acts freely on the manifold $M=\C^*=\C\setminus\{0\}$ 
by multiplication by the primitive $r$-th root of unity 
$\zeta=\exp(2\pi\ii/r)$.  In this instance, Theorem~\ref{thm:slfmono} 
provides a useful alternative to the above result.  This orbit 
configuration space is given by
\[
F_\Gamma(\C^*,n)=\{(x_1,\dots,x_n) \in (\C^*)^n \mid x_j \neq \zeta^p 
x_i 
\ \text{for $i\neq j$ and $1\le p \le r$}\},
\]
and thus may be realized as as the complement in $\C^n$ of the 
arrangement $\A_{r,n}$ consisting of the hyperplanes $H_j=\ker(x_j)$, 
$1\le j\le n$, and $H_{i,j}^{(p)}=\ker(x_j- \zeta^p x_i)$, $1\le i<j 
\le n$, $1\le p \le r$.  These are the reflecting hyperplanes of the 
full monomial group $G(r,n)$, the complex reflection group isomorphic 
to the wreath product of the symmetric group $\Sigma_n$ and 
$\Gamma=\Z/r\Z$, so we call $\A_{r,n}$ the (full) monomial 
arrangement.  The projection $p_\Gamma:F_\Gamma(\C^*,n)\to 
F_\Gamma(\C^*,n-1)$ reveals the fiber-type structure of this 
arrangement: $\A_{r,k+1}$ is strictly linearly fibered over $\A_{r,k}$ 
for $1\le k<n$.  For the arrangement $\A_{r,k+1}$, the root map 
$g_k:F_\Gamma(\C^*,k)\to F(\C,rk+1)$ of \eqref{eq:rootmap} is given by
\begin{equation} \label{eq:MonomialRootMap}
g_k(x_1,\dots,x_k)=(0,\z x_1,\dots,\z^r x_1,\z x_2,\dots,\z^r 
x_2,\dots\dots,
\z x_k,\dots,\z^r x_k).
\end{equation}
We explicitly record the results of Theorems~\ref{thm:slfmono} and 
\ref{thm:ftmono} in this special case.

\begin{thm} \label{thm:MonomialBundleThm}
For $\Gamma=\Z/r\Z$ acting freely on $\C^*$, the orbit configuration 
space bundle $p_\Gamma:F_\Gamma(\C^*,n+1) \to F_\Gamma(\C^*,n)$ is 
equivalent to the pullback of the classical configuration space bundle 
$p_{rn+2}:F(\C,rn+2)\to F(\C,rn+1)$ along the map $g_n$.  
Consequently,
\begin{enumerate}
\item the bundle $p_\Gamma:F_\Gamma(\C^*,n+1) \to F_\Gamma(\C^*,n)$ 
admits a section; 

\item the structure group of the bundle $p_\Gamma:F_\Gamma(\C^*,n+1) 
\to F_\Gamma(\C^*,n)$ is the Artin pure braid group on $rn+1$ strands 
$P_{rn+1}$;

\item the monodromy of $p_\Gamma:F_\Gamma(\C^*,n+1) \to 
F_\Gamma(\C^*,n)$ factors as $\a_{rn+1}\circ \c_n$, where 
$\c_n:\pi_1(F_\Gamma(\C^*,n)) \to P_{rn+1}$ is the map on fundamental 
groups induced by $g_n$, and $\a_{rn+1}:P_{rn+1} \to \Aut(\F_{rn+1})$ 
is the Artin representation;

\item the arrangement $\A_{r,n}$ is $K(\pi,1)$, and the group 
$\pi_1(F_\Gamma(\C^*,n)) = \rtimes_{j=1}^{n} \F_{r(j-1)+1}$ admits the 
structure of an iterated semidirect product of free groups.
\end{enumerate}
\end{thm}

\subsection{Pure Monomial Braid Groups} \label{subsec:puremono}
\ 
\medskip

In this section, we investigate the structure of the fundamental group 
of the orbit configuration space $F_\Gamma(\C^*,n)$, where 
$\Gamma=\Z/r\Z$.  Since this space is the complement of the reflection 
arrangement associated to the monomial group $G(r,n)$, we call this 
fundamental group the {\em pure monomial braid group}, and write 
$P(r,n)=\pi_1(F_\Gamma(\C^*,n))$.

We first construct some geometric braids in the orbit configuration 
space $F_\Gamma(\C^*,n)$.  A classical braid on $n$ strands may be 
described as (an equivalence class of) the motion of $n$ distinct 
points in the plane through time.  Thus a braid $\beta$ may be 
represented by a collection of $n$ maps 
$\beta(t)=(b_1(t),\dots,b_n(t))$, $b_j:[0,1]\to\C$, which satisfy 
$\beta(0)=\beta(1)$ (as sets), and $b_i(t)\neq b_j(t)$ for all $t$ if 
$i\neq j$.  If $b_j(0)=b_j(1)$ for each $j$, the braid $\beta$ is 
pure, and represents an element of the fundamental group of the 
configuration space $F(\C,n)$.  We adapt these ideas to the orbit 
configuration space $F_\Gamma(\C^*,n)$.

The root map $g_n:F_\Gamma(\C^*,n) \to F(\C,rn+1)$ of 
\eqref{eq:MonomialRootMap} is an imbedding, and the orbit 
configuration space $F_\Gamma(\C^*,n)$ is homeomorphic to the section 
$S(r,n)$ of $F(\C,rn+1)$ defined by the image of $g_n$.  If $\b$ is a 
braid in the configuration space $F(\C,rn+1)$, call $\b$ a {\em 
monomial braid} if $\b(t)\in S(r,n)$ for all $t$.  We give two 
relevant examples.

\begin{figure}[h]
\epsfysize=1 truein \epsfbox{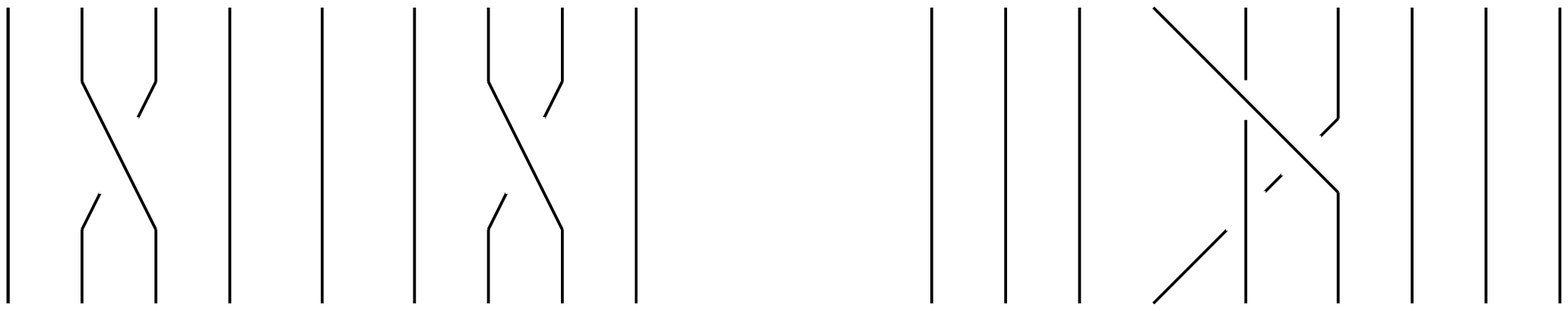}
\caption{Monomial braids $\rho_2$ (left) and $\rho_0$ (right) for 
$r=2$}
\label{fig:mb2}
\end{figure}

\begin{figure}[h]
\epsfysize=1.85 truein \epsfbox{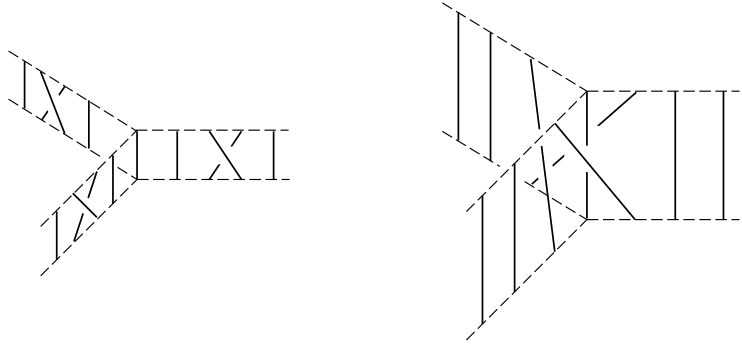}
\caption{Monomial braids $\rho_2$ (left) and $\rho_0$ (right) for 
$r=3$}
\label{fig:mb3}
\end{figure}

Define monomial braids $\rho_i=\rho_i(t)$ by
\begin{equation} \label{eq:FullMonomialGens}
\begin{split}
\rho_0&=g_n\circ (\exp(2\pi t\ii/r),2,3,\dots,n),\ \text{and}\\
\rho_i&=g_n\circ (1,\dots,i-1,i-\exp(\pi t\ii),i+1+\exp(\pi 
t\ii),i+2,\dots,n) 
\end{split}
\end{equation}
for $1\le i < n$.  Pictures of these braids, for small $r$, are given 
in Figures~\ref{fig:mb2} and \ref{fig:mb3}.

Expressing these monomial braids in terms of the standard generators 
$\s_i$ of the Artin braid group $B_{rn+1}$ is in general fairly 
difficult.  We describe one way to accomplish this.

For each $i$, the strands of the braid $\rho_i$ emanate from the 
points $0,k\z^p \in \C$, $1\le k \le n$, $1\le p \le r$.  Order these 
$rn+1$ strands as follows:
\begin{equation} \label{eq:strands}
\begin{matrix}
\text{strand \#}:\quad 1&2&3&\dots&r&r+1&r+2&\dots&2r+1&\dots\\
\hfill \text{point}:\quad 
0&\z&\z^2&\dots&\z^{r-1}&\zeta^r&2\z&\dots&2\z^r&\dots
\end{matrix}
\end{equation}
For $p=1,\dots,r-1$, successively rotate the rays $\z^{r-p} \cdot 
\R^+$ until the initial points of the strands lies on the positive 
real axis in order.  With these choices, the monomial braids are given 
by
\begin{equation} \label{eq:ArtinMonomialBraids}
\begin{split}
\rho_0 &= \s_r\s_{r-1}\cdots \s_1\s_1,\\
\rho_1 &= \tau_1^{-1}\cdot \s_2\s_4\cdots \s_{2r}\cdot \tau_1\\
\rho_{i} &= \tau_{i}^{-1}\cdot\s_{(i-1)r+2}\s_{(i-1)r+4}\cdots
\s_{ir+r} \cdot \tau_{i},\quad\text{for $2\le i\le n-1$,}
\end{split}
\end{equation}
where 
\[
\begin{align*}
\tau_1&=(\s_3\s_5\cdots\s_{2r-1})
(\s_4\s_6\cdots\s_{2r-2})
\cdots (\s_{r-1}\s_{r+1}\s_{r+3})
(\s_r\s_{r+2})(\s_{r+1}),\ \text{and for $i>1$,}\\ 
\tau_{i}&=(\s_{(i-1)r+3}\s_{(i-1)r+5}\cdots\s_{ir+r-1})
\cdots (\s_{ir-1}\s_{ir+1}\s_{ir+3})
(\s_{ir}\s_{ir+2})(\s_{ir+1}).
\end{align*}
\]
Using this, or experimenting with geometric braids, one can show 
that the monomial braids $\rho_0,\rho_1,\dots,\rho_{n-1}$ satisfy the
relations
\begin{equation} \label{eq:MonomialBraidRels}
\begin{split}
(\rho_0\rho_1)^2=(\rho_1\rho_0)^2,&\qquad  
\rho_i\rho_{i+1}\rho_i=\rho_{i+1}\rho_i\rho_{i+1}\quad (1\le 
i<n),\quad \text{and}\\  
&\rho_i\rho_j=\rho_j\rho_i\quad (|j-i|>1).
\end{split}
\end{equation}

Let $B(r,n)$ denote the group with generators 
$\rho_0,\rho_1,\dots,\rho_{n-1}$ and relations 
\eqref{eq:MonomialBraidRels}.  This is the (full) monomial braid 
group, the fundamental group of the quotient space 
$F_\Gamma(\C^*,n)/W$, where $W=G(r,n)$ is the full monomial group, 
cf.~\cite{BMR}.  Note that $B(r,n)$ is independent of $r$.  This group 
admits natural surjection to $G(r,n)$, which may be presented with 
generators $\rho_0,\rho_1,\dots,\rho_{n-1}$ and relations 
\eqref{eq:MonomialBraidRels} together with
\[
\rho_0^r=\rho_1^2=\dots =\rho_{n-1}^2=1.
\]

The pure monomial braid group may be realized as 
$P(r,n)=\ker(B(r,n)\to G(r,n))$, the kernel of the aforementioned 
surjection, see \cite{BMR}.  Elements of $P(r,n)$ are represented by 
(equivalence classes of) monomial braids $\beta$ as above with the 
property that $b_j(0)=b_j(1)$ for all $j$, thus are elements of the 
fundamental group of the space $S(r,n)$, homeomorphic to the orbit 
configuration space $F_\Gamma(\C^*,n)$.  We specify a number of these 
pure monomial braids.

For $1\le i \le n$, let 
$X_i=\rho_{i-1}\cdots\rho_2\rho_1\rho_0\rho_1\rho_2 \cdots\rho_{i-1}$, 
and define
\begin{equation} \label{eq:PureMonomialGens}
\begin{split} 
Z_j^{}&=\rho_{j-1}^{}\cdots\rho_2^{}\rho_1^{}\rho_0^r
\rho_1^{-1}\rho_2^{-1}\cdots\rho_{j-1}^{-1}
\ (1\le j\le n),\\
A_{i,j}^{(p)}&=X_i^{p-r} \cdot\rho_{j-1}^{}\cdots\rho_{i+1}^{} 
\rho_i^2
\rho_{i+1}^{-1}\cdots\rho_{j-1}^{-1}\cdot X_i^{r-p}
\ (1\le i<j\le n,\ 1\le p\le r).
\end{split}
\end{equation}

\begin{prop} \label{prop:PureMonomialGens}
The pure monomial braids $Z_j$ and $A_{i,j}^{(p)}$ $(1\le i<j,\ 1\le 
p\le r)$ generate the factor $\F_{r(j-1)+1}$ in the realization 
$P(r,n)=\rtimes_{j=1}^{n} \F_{r(j-1)+1}$ of the pure monomial braid 
group as an iterated semidirect product of free groups.  In 
particular, these braids (for $1\le j \le n$) generate the group 
$P(r,n)$.
\end{prop}
\begin{proof}
It suffices to show that the braids $Z_n$, $A_{i,n}^{(p)}$ generate 
the free group $\F_{r(n-1)+1}=\ker(P(r,n)\to P(r,n-1))$.  Identify the 
pure monomial braid group with the fundamental group of the space 
$S(r,n)$.

The free group $\F_{r(n-1)+1}$ is the fundamental group of the fiber 
of the projection $p_\Gamma:F_\Gamma(\C^*,n)\to F_\Gamma(\C^*,n-1)$.  
Via the imbeddings $g_k:F_\Gamma(\C^*,k)\to F(\C,rk+1)$, this 
projection corresponds to the map $S(r,n)\to S(r,n-1)$ defined by 
forgetting the last $r$ coordinates.  In terms of geometric braids, 
this map is given by forgetting the last, or outermost, $r$ strands.  
 From this and the definitions of the braids $Z_n$ and 
$A_{i,n}^{(p)}$, it is clear that these braids are elements of 
$\F_{r(n-1)+1}=\ker(P(r,n)\to P(r,n-1))$.

Recall from \eqref{eq:strands} that the strands of monomial braids are 
indexed by the points $0,k\zeta^p \in \C$ from which these braids 
emanate.  Checking that the braid $Z_n$ performs a full twist on the 
$r$ strands emanating from $n\z,n\z^2,\dots,n\z^r$, and that 
$A_{i,n}^{(p)}$ simultaneously performs full twists on the pairs of 
strands emanating from $i\z^{p+q}$ and $n\z^q$, $1\le q\le r$, we see 
that the homology classes of these braids are independent.  It follows 
that the braids $Z_n$ and $A_{i,n}^{(p)}$ generate the free group 
$\F_{r(n-1)+1}$.
\end{proof}

In principle, Theorem \ref{thm:MonomialBundleThm} and the above result 
yield a presentation of the group $P(r,n)$ as an iterated semidirect 
product of free groups, see Remarks~\ref{rem:slfpres} and 
\ref{rem:ftpres}.  See Theorem \ref{thm:PBn} for the special case 
$r=2$.  In general, the requisite Artin representation calculations 
are too complicated to pursue fruitfully.  However, using Proposition 
\ref{prop:bmartin}, we can record some useful qualitative information 
concerning this presentation.  Recall that $\widetilde U$ denotes some 
conjugate of an element $U$ of a free group.

\begin{prop} \label{prop:RoughPureMonomialRels}
In the pure monomial braid group $P(r,n)=\rtimes_{j=1}^{n} 
\F_{r(j-1)+1}$, for $j<l$, the action of $\F_{r(j-1)+1}$ on 
$\F_{r(l-1)+1}$ is of the form
\[
\begin{alignat*}{2}
Z_j^{-1} Z_l^{} Z_j^{} &= U^{} \widetilde{Z}_l^{} U^{-1}, &\qquad
(A_{i,j}^{(p)})^{-1} Z_l^{} A_{i,j}^{(p)} &= \widetilde{Z}_l^{}, \\
Z_j^{-1} A_{k,l}^{(q)} Z_j^{} &= V_k^{}\widetilde{A}_{k,l}^{(q)} 
V_k^{-1},&\qquad
(A_{i,j}^{(p)})^{-1} A_{k,l}^{(q)} A_{i,j}^{(p)} &= 
W_k^{} \widetilde{A}_{k,l}^{(q)} W_k^{-1}, 
\end{alignat*}
\]
where $U=V_j=\widetilde{A}_{j,l}^{(1)} \widetilde{A}_{j,l}^{(2)} 
\cdots 
\widetilde{A}_{j,l}^{(r-1)} \widetilde{Z}_{l}^{} 
\widetilde{A}_{j,l}^{(r)}$,  $V_k=1$ if $k\neq j$, 
$W_i=\widetilde{A}_{i,l}^{(q)}\widetilde{A}_{j,l}^{(m)}$ if $q\equiv 
p+m \mod r$, 
$W_j=\widetilde{A}_{i,l}^{(m)}\widetilde{A}_{j,l}^{(q)}$ if $m\equiv 
p+q \mod r$, 
and $W_k=1$ otherwise.
\end{prop}

While the conjugates $\widetilde{Z}_{l}$ and 
$\widetilde{A}_{i,l}^{(q)}$ in the free group $\F_{r(l-1)+1}$ which 
appear above are not readily accessible, one can use this result and 
some general facts from the theory of hyperplane arrangements to 
obtain an explicit presentation of the group $P(r,n)$.  We briefly 
recall these facts concerning the ``Randell-Arvola'' presentation of 
an arrangement group.  See \cite{Ra, Ar, OT, CSbm} for detailed 
discussions.

For a general arrangement $\A=\{H_1,\dots,H_n\}$, the fundamental 
group $G=G(\A)$ of the complement is generated by meridonal loops, 
$g_j$ about $H_j$.  Relations in the group $G$ arise from codimension 
two intersections of hyperplanes.  If $H_1\cap \dots \cap H_m$ is such 
an intersection, then in $G$ one has a corresponding family of $m-1$ 
relations
\[
\widetilde{g}_1 \widetilde{g}_2\cdots \widetilde{g}_m = 
\widetilde{g}_2\cdots \widetilde{g}_m \cdot \widetilde{g}_1 = 
\cdots\cdots =
\widetilde{g}_m\cdot \widetilde{g}_1 \cdots \widetilde{g}_{m-1}.\]
This family of (commutation) relations is denoted by 
$[\widetilde{g}_1,\dots,\widetilde{g}_m]$.

For the monomial arrangement $A_{r,n}$, the codimension two 
intersections of hyperplanes are recorded implicitly in Proposition 
\ref{prop:RoughPureMonomialRels}.  They are, for $j<l$,
\[
\begin{alignat*}{2}
&H_j^{} \cap H_{j,l}^{(1)} \cap \dots \cap H_{j,l}^{(r-1)} \cap 
H_l^{} \cap
H_{j,l}^{(r)}, &\qquad 
&H_j^{}  \cap H_{k,l}^{(q)}\ \text{if}\ j\neq k,  
\qquad\qquad H_{i,j}^{(p)} \cap H_l^{},\\  
&H_{i,j}^{(p)} \cap H_{k,l}^{(q)}\ \text{if}\ 
\{i,j\}\cap\{k,l\}=\emptyset, &\qquad 
&H_{i,j}^{(p)} \cap H_{i,l}^{(m)} \cap H_{j,l}^{(q)} 
\ \text{if}\ m\equiv p+q \mod r.
\end{alignat*}
\]
While we do not record the conjugations arising in the iterated 
semidirect product structure of the group $P(r,n)$, we can, with some 
effort, record those arising in a Randell-Arvola presentation of this 
group.  Recall the monomial braids 
$X_i=\rho_{i-1}\cdots\rho_2\rho_1\rho_0\rho_1\rho_2 \cdots\rho_{i-1}$.  
For $p<r$, write $A_{i,j}^{[p]}=A_{i,j}^{(p)}A_{i,j}^{(p+1)} \cdots 
A_{i,j}^{(r-1)}$.  Set $A_{i,j}^{[r]}=1$.  We require the following 
technical result.  We omit the proof, which is a delicate exercise 
using the relations \eqref{eq:MonomialBraidRels} satisfied by monomial 
braids and the definition \eqref{eq:PureMonomialGens}.

\begin{lem} \label{lem:conj}
We have
\[
\begin{align*}
X_i^{-1} Z_{l} X_i^{} &=
\begin{cases}
Z_l&\text{if $l<i$,}\\
C_l^{-1} Z_l^{} C_l^{}&\text{if $l=i$,}\\
A_{i,l}^{(r-1)}Z_l(A_{i,l}^{(r-1)})^{-1}&\text{if $l>i$,}
\end{cases}\\
\intertext{
where $C_l=A_{1,l}^{(r)}A_{2,l}^{(r)}\cdots A_{l-1,l}^{(r)}$, and
}
X_i^{-1} A_{k,l}^{(r)} X_i^{} &=
\begin{cases}
A_{k,l}^{(r)}&\text{if $i<k$ or $i>l$,}\\
A_{k,l}^{(r-1)}&\text{if $i=k$,}\\
A_{i,l}^{(r-1)}A_{k,l}^{(r)}(A_{i,l}^{(r-1)})^{-1}&\text{if 
$k<i<l$,}\\
D_k^{}A_{k,l}^{(1)}D_k^{-1}&\text{if $i=l$,}
\end{cases}
\end{align*}
\]
where $D_k^{}=A_{k-1,k}^{[1]}A_{k-2,k}^{[1]}\cdots A_{1,k}^{[1]}Z_k^{}
C_k^{}$.
\end{lem}

Define elements $Q_{j,p} \in P(r,n)$ by $Q_{1,p}=1$, and for $j \ge 
2$,
\[
Q_{j,p}=C_j \prod_{q=1}^p \bigl( D_1^{} A_{1,j}^{(q)} D_1^{-1}
D_2^{} A_{2,j}^{(q)} D_2^{-1} \cdots D_{j-1}^{} A_{j-1,j}^{(q)} 
D_{j-1}^{-1}
\bigr).
\]

\begin{thm} \label{thm:PureMonomialPresentation}
The pure monomial braid group $P(r,n)$ admits a presentation with 
generators $Z_j$ $(1\le j \le n)$, $A_{i,j}^{(p)}$ $(1\le i<j\le n$, 
$1\le p\le r)$, and for $i<j<k<l$, relations
\begin{equation*}
\begin{align}
[Z_j^{},Q_{j,r-1}^{}A_{j,l}^{(1)}Q_{j,r-1}^{-1},\dots,
Q_{j,1}^{}A_{j,l}^{(r-1)}Q_{j,1}^{-1},Z_l^{},A_{j,l}^{(r)}],&
\label{eq:ZAZA}\\
[Z_i,A_{k,l}^{(p)}],\quad [A_{i,j}^{(p)},Z_k], \quad
[A_{i,j}^{[p]}Z_j(A_{i,j}^{[p]})^{-1},
(A_{i,k}^{(p)})^{(A_{i+1,k}^{(r)}\cdots A_{j,k}^{(r)})}],
&\quad\text{$1\le p \le r$,}\label{eq:ZAfirst}\\
[A_{i,j}^{(p)},A_{k,l}^{(q)}],\quad 
[A_{j,k}^{(p)},A_{j,l}^{[p]}A_{i,l}^{(q)}(A_{j,l}^{[p]})^{-1}],
&\quad\text{$1\le p,q\le r$,}\label{eq:AAfirst}\\
[A_{j,k}^{[p]}A_{i,k}^{(q)}(A_{j,k}^{[p]})^{-1},
(A_{j,l}^{(p)})^{A_{k,l}^{(r)}}], \quad 
[A_{i,j}^{(p)},A_{i,k}^{(p)},A_{j,k}^{(r)}],
&\quad\text{$1\le p,q\le r$,}\\
[A_{i,j}^{(p)},A_{j,k}^{(r-q)},
A_{j,k}^{[r-q+1]}A_{i,k}^{(p-q)}
(A_{j,k}^{[r-q+1]})^{-1}],
&\quad\text{$1\le q<p\le r$,}\\
[D_i^{}A_{i,j}^{(p)}D_i^{-1},A_{j,k}^{(r-q)},
A_{j,k}^{[r-q+1]}A_{i,k}^{(r-q+p)}
(A_{j,k}^{[r-q+1]})^{-1}],
&\quad\text{$1\le p \le q< r$,}\label{eq:AAlast}
\end{align}
\end{equation*}
\end{thm}
\begin{proof}[Sketch of Proof]
It follows from the above discussion that the group $P(r,n)$ admits a 
presentation of this form.  We sketch how one may use Lemma 
\ref{lem:conj} and the monomial braid relations 
\eqref{eq:MonomialBraidRels} to show that the conjugating words 
appearing in the relation families \eqref{eq:ZAZA}--\eqref{eq:AAlast} 
are as asserted.

Using \eqref{eq:MonomialBraidRels}, one can show that 
$Z_1^{}A_{1,2}^{(1)}\cdots A_{1,2}^{(r-1)}Z_2^{}A_{1,2}^{(r)}= 
(\rho_0\rho_1)^{2r}$ commutes with $Z_1$, $Z_2$, and $A_{1,2}^{(p)}$ 
for each $p$.  Relations \eqref{eq:ZAZA} with $(j,l)=(1,2)$ follow.  
Let $\nu_2=\rho_1\rho_2$ and $\nu_j=\rho_1\cdots\rho_j\cdot\nu_{j-1}$ 
for $j>2$.  Then $\nu_j^{} Z_1^{} \nu_j^{-1}=Z_j^{}$ and $\nu_j Z_2^{} 
\nu_j^{-1}=Z_{j+1}^{}$.  Conjugating \eqref{eq:ZAZA} with 
$(j,l)=(1,2)$ by $\nu_j$ yields \eqref{eq:ZAZA} for any $j$ and 
$l=j+1$.  Conjugating these last relations by $\rho_l\cdots\rho_{j+1}$ 
yields \eqref{eq:ZAZA} in general.

Similarly, using \eqref{eq:MonomialBraidRels}, one can show that the 
following relations hold in $P(r,n)$:
\[
[Z_i,A_{j,k}^{(r)}],\quad [Z_j,
(A_{i,k}^{(r)})^{(A_{i+1,k}^{(r)}\cdots A_{j,k}^{(r)})}],\quad 
\text{and} \quad
[A_{i,j}^{(r)},Z_k].
\]
Repeated conjugation these relations by $X_t$, $t=i,j,k$, and use of 
Lemma \ref{lem:conj} yields the relations \eqref{eq:ZAfirst}.  

Finally, it is clear from the definition \eqref{eq:PureMonomialGens} 
that the pure monomial braids $A_{i,j}^{(r)}$ satisfy the classical 
pure braid relations \eqref{eq:PureBraidRels}.  These relations may 
be rewritten as
\[
[A_{i,j}^{(r)},A_{k,l}^{(r)}],\quad[A_{j,k}^{(r)},A_{i,l}^{(r)}],
\quad [A_{i,k}^{(r)},(A_{j,l}^{(r)})^{A_{k,l}^{(r)}}], \quad 
\text{and} \quad [A_{i,j}^{(r)},A_{i,k}^{(r)},A_{j,k}^{(r)}].
\] 
Repeated conjugation these relations by $X_t$, $t=i,j,k,l$, and use 
of Lemma \ref{lem:conj} yields the relations 
\eqref{eq:AAfirst}--\eqref{eq:AAlast}.
\end{proof}

\begin{rem}
In the special case $r=2$, one can check that the presentations of the 
generalized pure braid group $PB_n=P(2,n)$ given in Theorems 
\ref{thm:PBn} and \ref{thm:PureMonomialPresentation} are equivalent.  
The correspondence between the generators in these two results is 
given by $a_{i,j}=A_{i,j}^{(2)}$, $b_{i,j}=A_{i,j}^{(1)}$, and 
$c_j=Z_j$.
\end{rem}

\subsection{The Lie Algebra Associated to the Lower Central Series} 
\label{subsec:LCS}
\ 
\medskip

We now show how the results of the previous section may be used to 
determine the structure of the Lie algebra associated to the lower 
central series of the pure monomial braid group.

For any group $G$, let $G_k$ denote the $k$-th lower central series 
subgroup, defined inductively by $G_1=G$ and $G_{k+1}=[G_k,G]$ for $k 
\ge 1$.  Let $\gr G = \bigoplus_{k \ge 1} G(k)$, where 
$G(k)=G_k/G_{k+1}$.  The map (of sets) $G \times G \to G$ given by the 
commutator, $(x,y) \mapsto [x,y]=xyx^{-1}y^{-1}$, induces a bilinear 
map $\gr G \times \gr G \to \gr G$ which defines a Lie algebra 
structure on $\gr G$, see for instance \cite[Chapter I]{Serre}.  

\begin{exm} \label{exm:free}
If $G=\F_n$ is a finitely generated free group, then $\gr G$ is 
isomorphic to the free Lie algebra $L(n)$ on $n$ generators, see 
\cite[Chapter IV]{Serre}.
\end{exm}

The additive structure of the Lie algebra associated to the lower 
central series of the fundamental group of the complement of a 
fiber-type arrangement is given by the following result.

\begin{thm}[Falk and Randell \cite{FR1}, Theorem~3.1] \label{thm:FR}
Let $1 \to H 
\xrightarrow{i} G \xrightarrow{j} K \to 1$ be a split extension of 
groups such that $K$ acts trivially on $H_1/H_2$.  Then the induced 
sequence of graded abelian groups, $0 \to \gr H \to \gr G \to 
\gr K \to 0$, is split exact.
\end{thm}

In the case where $G$ is the group of a fiber-type arrangement, the 
proof of this result in \cite{FR1} shows that the lower central series 
quotients $G(k)$ are free abelian for all $k$, see \cite{FR2}.  
Applying this result and observation inductively, we obtain

\begin{cor} \label{cor:grFT}
Let $\A$ be a fiber-type arrangement with exponents 
$\{d_1,\dots,d_\ll\}$ and group $G$.  Then, as abelian groups, $\gr G 
\cong L(d_1) \oplus L(d_2) \oplus \dots \oplus L(d_\ll)$.
\end{cor}

We now pursue the Lie bracket relations in $\gr G$ for these groups, 
in particular for the pure monomial braid groups $P(r,n)$.

In the situation of Theorem \ref{thm:FR}, the group $G \cong H 
\rtimes_\a K$ may be realized as the semidirect product of $H$ and 
$K$, determined by a homomorphism $\a:K \to \Out(H)$.  If $s:K\to G$ 
denotes the splitting, then $G$ is generated by $i(H)$ and $s(K)$.  
Identifying $H$ and $K$ with their images in $G$, we have relations 
$x^{-1}y^{}x^{}=\a(x)(y)$ in $G$ for $y\in H$ and $x\in K$.  If the 
action of $K$ on $H$ is by conjugation, the Lie bracket relations in 
$\gr G$ are readily obtained from the relations in $G$ itself, as 
follows.

\begin{lem} \label{lem:brackets}
Let $y\in H$ and $x\in K$.  If $x^{-1}y^{}x^{}=w^{}y^{}w^{-1}$ for 
some $w\in H$, then $[\bar{x}+\bar{w},\bar{y}]=0$ in $\gr G$, where 
$\bar{x}$, $\bar{y}$, and $\bar{w}$ are the images of $x$, $y$, and 
$w$ in $\gr G$.
\end{lem}
\begin{proof}
Rewriting $x^{-1}y^{}x^{}=w^{}y^{}w^{-1}$ as 
$1_G=x^{}w^{}y^{}w^{-1}x^{-1}y^{-1}=[xw,y] \in G_2$, the result 
follows from the commutator identity 
$\bigl[xw,y\bigr]=\bigl[x,[w,y]\bigr]\bigl[w,y\bigr]\bigl[x,y\bigr]$ 
in $G$.
\end{proof}

We subsequently write simply $g$, as opposed to $\bar{g}$, for the 
image in $\gr G$ of $g\in G$.

\begin{exm}
The structure of the Lie algebra associated to the lower central 
series of the classical pure braid group $P_n=\rtimes_{j=1}^{n-1} 
\F_j$ was determined by Kohno \cite{Ko}.  The above considerations may 
be used to recover this result.  By Corollary \ref{cor:grFT}, we have 
$\gr P_n \cong \bigoplus_{j=1}^{n-1} L(j)$, where the free Lie algebra 
$L(j)$ is generated by $A_{1,j+1},\dots,A_{j,j+1}$.  Applying Lemma 
\ref{lem:brackets} to the defining relations \eqref{eq:PureBraidRels} 
of $P_n$, we see that the Lie bracket relations in $\gr P_n$ are the 
``infinitesimal pure braid relations,'' given by
\[
\begin{align*}
[A_{i,j}+A_{i,k}+A_{j,k},A_{m,k}]&=0 \quad
\text{for $m=i,j$, and}\\
[A_{i,j},A_{k,l}]&=0
\quad\text{for $\{i,j\}\cap\{k,l\}=\emptyset$.}
\end{align*}
\]

The Lie algebra $\gr P_n$ arises in a number of other contexts.  For 
instance, the integral homology, $H_*(\Omega F(\C^k,n))$, of the loop 
space of the classical configuration space was recently computed by 
Fadell and Husseini \cite{FH}.  Subsequently, Cohen and Gitler 
\cite{CG} showed that, with appropriate regrading, the Lie algebra of 
primitives, $PH_*(\Omega F(\C^k,n))$, is isomorphic to $\gr P_n$ for 
$k\ge 2$.
\end{exm}

The orbit configuration spaces $F_\Gamma(\C^k \setminus \{0\},n)$, 
where $\Gamma=\Z/r\Z$ acts on $\C^k \setminus \{0\}$ by multiplication 
by a primitive $r$-th root of unity, provide natural generalizations 
of these results.  In the case $k=1$, using 
Theorem~\ref{thm:PureMonomialPresentation} (or 
Proposition~\ref{prop:RoughPureMonomialRels}), we obtain:

\begin{thm} \label{thm:PureMonomialBrackets}
Let $\gr P(r,n)$ be the Lie algebra associated to the lower central 
series of the pure monomial braid group.  Then, $\gr P(r,n) \cong 
\bigoplus_{j=0}^{n-1} L(rj+1)$ as abelian groups, where $L(rj+1)$ is 
generated by $Z_{j+1}$ and $A_{i,j+1}^{(p)}$, $1\le i\le j$, $1\le p 
\le r$.  The Lie bracket relations in $\gr P(r,n)$ are given by
\[
\begin{align*}
\bigl[Z_j+Z_l+A_{j,l}^{(1)}+A_{j,l}^{(2)}+\dots+A_{j,l}^{(r)},Y\bigr]
&=0 
\quad
\text{for $Y=Z_l$, $Y=A_{j,l}^{(p)}$, $1\le p\le r$,}\\
[A_{i,j}^{(p)}+A_{i,k}^{(q)}+A_{j,k}^{(m)},Y]&=0
\quad\text{for $Y=A_{i,k}^{(q)},A_{j,k}^{(m)}$, $q\equiv p+m \mod 
r$,}\\
[A_{i,j}^{(p)},A_{k,l}^{(q)}]&=0
\quad\text{for $\{i,j\}\cap\{k,l\}=\emptyset$, $1\le p,q\le r$, and}\\
[Z_{k},A_{i,j}^{(p)}]&=0
\quad\text{for $k\neq i,j$ and $1\le p \le r$.}
\end{align*}
\]
\end{thm}

For $k\ge 2$, consider the loop space $\Omega 
F_\Gamma(\C^k\setminus\{0\},n)$.  The Lie algebra of primitives in the 
homology, $H_*(\Omega F_\Gamma(\C^k\setminus\{0\},n))$, was calculated 
by Xicot\'encatl \cite[Theorem 3.1.2]{Xi}.  Comparing the Lie bracket 
relations obtained there and those in $\gr P(r,n)$ recorded above, we 
obtain the following, which was conjectured in \cite[Section 3.7]{Xi}.

\begin{cor} \label{cor:Xconj}
For $k\ge 2$, with appropriate regrading, the Lie algebra $\gr P(r,n)$ 
is isomorphic to the Lie algebra of primitives, $PH_*(\Omega 
F_\Gamma(\C^k\setminus\{0\},n))$, in the Hopf algebra $H_*(\Omega 
F_\Gamma(\C^k\setminus\{0\},n))$.
\end{cor}

Similar considerations reveal the structure of the Lie algebra 
associated to the lower central series of the fundamental group of the 
complement of an arbitrary fiber-type arrangement $\A$ defined by 
$Q(\A)=\prod_{j=1}^{\ll} Q_{j}(x_{1},\dots,x_{j})$.  Recall the 
presentation of the group $G$ of such an arrangement from 
Remark~\ref{rem:ftpres}, with generators $x_{q,j}$ and relations $ 
x_{p,i}^{-1}x_{q,j}^{}x_{p,i}^{}=\eta_j^{}(x_{p,i}^{})(x_{q,j}^{})$.  
The generators $x_{q,j}$ correspond to the hyperplanes $H_{q,j}$ 
defined by the degree $d_{j}$ polynomial $Q_{j}(x_{1},\dots,x_{j})$, 
and $\eta_j=\a_{d_j} \circ \c_j$ is the composition of the Artin 
representation and the homomorphism induced by the map $g_j$ 
of~\eqref{eq:rootmap}.

 From the identification of the monodromy of the strictly linearly 
fibered bundle $M(\A_{j}) \to M(\A_{j-1})$ in Section 
\ref{subsec:monocalc}, for $i<j$ we have $\c_j(x_{p,i})=\prod_{k=1}^r 
A_{V_k}^{\zeta_k}$, where $V=(V_{1}\mid V_{2}\mid \dots \mid V_{r})$ 
is the partition of $[d_j]$ recording the hyperplanes of 
$\A_j\setminus\A_{j-1}$ which meet $H_{p,i}$ in codimension two, 
cf.~Theorem \ref{thm:bmpure} and Remark \ref{rem:mono&comb}.  Applying 
the Artin representation as in Proposition \ref{prop:bmartin}, the 
relations in $G$ may be expressed more explicitly as
\begin{equation} \label{eq:ftrels}
x_{p,i}^{-1}\cdot t_{q}^{}\cdot x_{p,i}^{}=
\widetilde{t}_{V_k}^{} \cdot \widetilde{t}_{q}^{} \cdot 
\widetilde{t}_{V_k}^{-1},
\end{equation}
where $t_m=x_{m,j}$ for $1\le m \le d_j$ and $V_k$ is the unique block 
of $V$ containing $q$.

\begin{thm} \label{thm:grftG}
Let $\A$ be a fiber-type arrangement with exponents 
$\{d_1,\dots,d_\ll\}$, and group $G=G(\A)$.  Then, $\gr G \cong 
\bigoplus_{j=1}^\ll L(d_j)$ as abelian groups, where $L(d_j)$ is 
generated by $x_{q,j}$, $1\le p \le d_j$, $1\le j\le \ll$.  The Lie 
bracket relations in $\gr G$ are given by
\begin{equation}\label{eq:ftbrackets} 
[x_{p,i}+x_{q_1,j}+x_{q_2,j}+ 
\dots+x_{q_m,j},x_{q_k,j}]=0 \quad \text{for}\ 1\le k \le m,
\end{equation}
where $1\le i<j \le \ll$ and $\{q_1,\dots,q_m\}$ is maximal such that 
$\codim H_{p,i} \cap \bigcap_{k=1}^m H_{q_k,j} = 2$.
\end{thm}
\begin{proof}
The isomorphism, $\gr G \cong \bigoplus_{j=1}^\ll L(d_j)$, of graded 
abelian groups was noted in Corollary \ref{cor:grFT} above.  The Lie 
bracket relations in $\gr G$ may be obtained by applying Lemma 
\ref{lem:brackets} to the relations \eqref{eq:ftrels}.
\end{proof}

For any arrangement $\A$ of $n$ hyperplanes, the (rational) holonomy 
Lie algebra $\LL_{\Q}$ of the complement $M$ is the quotient of the 
free Lie algebra $L_{\Q}(n)=L(H_{1}(M;\Q))$ by the image of the map 
$H_{1}(M;\Q) \to \bigwedge^{2}H_{1}(M;\Q)$ dual to the cup product.  
For a fiber-type arrangement with exponents $\{d_{1},\dots,d_{\ll}\}$, 
Jambu \cite{Ja} shows that $\LL_{\Q} \cong \bigoplus_{j=1}^{\ll} 
L_{\Q}(d_{j})$ as graded vector spaces.  For an arbitrary arrangement, 
Kohno \cite{K1} shows that $\LL_{\Q}$ is generated by elements 
$x_{1},\dots,x_{n}$ in one-to-one correspondence with the hyperplanes 
of $\A$, with relations
\begin{equation} \label{eq:holonomybrackets}
[x_{q_{1}}+x_{q_{2}}+\dots+x_{q_{m}},x_{q_{k}}]=0, \quad 1\le k \le m,
\end{equation}
for each maximal family $\{H_{q_{1}},\dots,H_{q_{m}}\}$ of hyperplanes 
of $\A$ with $\codim \bigcap_{k=1}^{m} H_{q_{k}} = 2$.

In light of this result, by imposing the relations 
\eqref{eq:holonomybrackets} on the free Lie algebra 
$L(n)=L_{\Z}(n)=L(H_{1}(M;\Z))$, we may consider the integral holonomy 
Lie algebra $\LL=\LL_{\Z}$ of $M$.  Furthermore, comparing the 
relations \eqref{eq:ftbrackets} and \eqref{eq:holonomybrackets} for a 
fiber-type arrangement, we have the following result, which may also 
be obtained from the work of Kohno \cite{K1,Ko}, together with the 
results of Falk and Randell stated above.

\begin{cor}
Let $\A$ be a fiber-type arrangement with complement $M$ and group 
$G$.  Then the integral holonomy Lie algebra $\LL$ of $M$ is 
isomorphic to the Lie algebra $\gr G$ associated to the lower central 
series of $G$.
\end{cor}

\begin{ack}
Portions of this work were carried during visits to the Department of 
Mathematics at the University of Wisconsin-Madison in the Spring of 
1999.  We thank the Department for its hospitality, and the College of 
Arts {\&} Sciences at Louisiana State University for granting the 
Research Fellowship which made these visits possible.  We also thank 
Fred Cohen and Peter Orlik for useful conversations.  Some of the 
results presented here stem from our collaboration with Alex Suciu, 
who we acknowledge.
\end{ack}

\bibliographystyle{amsalpha}

\end{document}